\documentclass{amsart}
\usepackage{amssymb}
\usepackage{amsfonts}
\usepackage{amstext}

\numberwithin{equation}{section}
\newtheorem{theorem}{Theorem}[section]
\newtheorem{lemma}[theorem]{Lemma}
\newtheorem{proposition}[theorem]{Proposition}
\newtheorem{corollary}[theorem]{Corollary}

\theoremstyle{definition}
\newtheorem{definition}[theorem]{Definition}
\newtheorem{example}[theorem]{Example}

\theoremstyle{remark}
\newtheorem{remark}[theorem]{\bf{Remark}}

\newcommand{\R}{{\mathbb{R}}}

\newcommand{\Z}{{\mathbb{Z}}}

\newcommand{\CI}{{\mathcal{I}}}
\newcommand{\CL}{{\mathcal{L}}}
\newcommand{\CD}{{\mathcal{D}}}

\newcommand{\<}{{\langle}}
\renewcommand{\>}{{\rangle}}

\newcommand{\tens}{\otimes}
\newcommand{\id}{{\rm id}}
\newcommand{\extd}{{\rm d}}

\newcommand{\la}{{\triangleright}}

\def\lcross{{>\!\!\!\triangleleft}}

\renewcommand{\o}{{}_{\scriptscriptstyle(1)}}
\renewcommand{\t}{{}_{\scriptscriptstyle(2)}}

\newcommand{\uo}{{}^{\scriptscriptstyle 1}}
\newcommand{\ut}{{}^{\scriptscriptstyle 2}}

\newcommand{\lo}{{[\![}}
\newcommand{\lc}{{]\!]}}

\newcommand{\inter}{{\mathfrak{i}}}
\newcommand{\rinter}{{\grave{\mathfrak{i}}}}
\newcommand{\cj}{{\mathfrak{j}}}

\begin{document}

\title{Reconstruction and quantization of Riemannian structures}


\author{S. Majid}
\address{Queen Mary University of London\\
School of Mathematics, Mile End Rd, London E1 4NS, UK}

\email{s.majid@qmul.ac.uk}

\date{December 2013, Ver 3}


\begin{abstract} We study how the Riemannian structure on a manifold can be usefully reconstructed from its codifferential $\delta$, including a formula  $\nabla_\omega\eta={1\over 2}( \delta(\omega\eta)-(\delta\omega)\eta+\omega(\delta\eta)  +\CL_\omega\eta+\inter_\eta\extd\omega)$ for the Levi-Civita covariant derivative in terms of 1-forms, where $\CL, \inter$ are respectively the Lie derivative and interior product along the corresponding vector fields. The covariant derivative extends naturally along forms of any degree and to possibly degenerate $(\ ,\ )$. In the nondegenerate case, $\delta$ makes the exterior algebra into a BV algebra.  In the invertible case we show that ${\rm Ricci}=-{1\over 2}\Delta g$ where the Hodge Laplacian $\Delta$ extends in a natural way to act on the metric. 

Our results come from a new way of thinking about metrics and connections as a kind of cocycle data for central extensions of differential graded algebras (DGAs), a theory which we introduce. We  show that any cleft extension of the exterior algebra  $\Omega(M)$ on a manifold is associated to a possibly-degenerate metric and covariant derivative. Those for which $\extd$ is not deformed up to isomorphism correspond to the Levi-Civita case. We provide a construction for such extensions both of classical DGAs and of already non-graded-commutative DGAs, thereby constructing a class of bimodule covariant derivatives via a kind of quantum analogue of the Koszul formula. 

We also provide a semidirect product of any differential graded algebra by the quantum differential algebra $\Omega(t,\extd t)$ in one variable, to introduce a noncommutative `time'. Composing these two constructions recovers a differential quantisation of $M\times\R$ in  \cite{Ma:bh}.  \end{abstract}
\maketitle

\section{Introduction}

We describe a new way of thinking about Riemannian geometry as coming out of algebra, at both the classical and the `quantum' or noncommutative geometry level.

Most of the paper is at the classical level of a manifold $M$ but for clarity done in terms of the graded-commutative algebra. Specifically, in Section~2 we address and answer affirmatively the following question: can a Riemannian structure on a manifold $M$ be usefully reconstructed from the algebraic properties of the codifferential $\delta$  on the exterior algebra $\Omega(M)$ of differential forms?  That one can do the reconstruction is not in doubt: the failure of $\delta$ to commute with functions easily allows one to recover the metric (Definition~\ref{regular}) after which the Levi-Civita covariant derivative is of course determined by the Koszul formula. What is not so obvious and which we find is that the Levi-Civita covariant derivative and its properties have a direct expression in these terms much in the style of the Cartan formula for the Lie derivative. Working with forms rather than vector fields (converting one to the other via the metric where needed) we find (Theorem~\ref{levi})
\[ \nabla_\omega\eta={1\over 2}(\delta(\omega\eta)-(\delta\omega)\eta+\omega\delta\eta+ \CL_\omega\eta +\extd\omega\perp\eta),\quad\forall \omega,\eta\in \Omega^1(M)\qquad\qquad (*)\]
which we explain as follows. First of all, it is well-known that the codifferential is not a graded-derivation of the exterior algebra but we see  that this failure equates to the antisymmetric part  of the Levi-Civita connection under interchange of $\omega,\eta$. The remaining terms are the symmetric part and these also have a natural expression in terms of the Lie derivative $\CL_\omega=\inter_\omega\extd+\extd\inter_\omega$ and an operation $\perp$ which in the present instance is the interior product $\inter$ but which more generally denotes its extension as a derivation in both arguments (see Section~2.1). The interior product itself is the metric inner product $(\ ,\ )$ on 1-forms extended as a graded-derivation. 

As well as providing an elegant formula, which we have not found elsewhere, our point of view leads to some immediate generalisations:

1) The formula (*) makes sense for $\eta$ of any degree and provides in fact the Levi-Civita covariant derivative on all degrees of $\Omega(M)$ at once.  

2) The formula (*) assumes a metric implicitly in the construction of $\delta$ but beyond that requires only the `metric inner product' $(\ ,\ )$ rather than the metric $g\in\Omega^1(M)\tens_{C^\infty(M)}\Omega^1(M)$ itself, providing a mild generalisation of classical Riemannian geometry. 

3) The formula (*) makes sense for covariant derivation `along' $\omega$ of any degree! This one we do not explore very much beyond noting some covariant derivative-like properties (Corollary~\ref{higher}) but it is this extended object that arises out of our algebraic point of view. As least in the nondegenerate case the extended $\nabla_\omega$ is a derivation of degree $|\omega|-1$. The Lie derivative along higher forms is defined by a higher-form Cartan formula and $\perp$ (see Theorem~\ref{cleftclasslevi}).

On the technical side it is known cf. \cite{Coll} that the failure of $\delta$ to be a graded-derivation (we will call such expressions `Leibnizators') is closely related to the Schouten bracket of alternating multivector fields \cite{Sch}. Unwinding this in degree 1,  (*) becomes closer to the Koszul formula as should be expected, and conversely the latter  in principle generalises to all degrees via the Schouten bracket. When $(\ ,\ )$ is not invertible we do not have such a conversion but if $(\ ,\ )$ is at least nondegenerate we find similar properties such as the Leibnizator of $\delta$ being a graded-derivation, which amounts to a  7-term triple product identity (Corollary~\ref{delta triple}) that makes $(\Omega(M),\delta)$ into BV algebra or a slight generalisation of it where $\delta^2$ is tensorial. We also solve the reconstruction problem as posed, finding (Proposition~\ref{bijection}) a bijection between codifferentials modulo interior codifferentials  and Riemannian structures. We use the divergence operator to show that every Riemannian structure has a representative $\delta$ with $\delta^2=0$, at least in the usual invertible case (Proposition~\ref{div}). 

4) As a modest application of our codifferential approach to Riemannian geometry we show that when $(\ ,\ )$ has an inverse, $g$, 
\[ {\rm Ricci}=-{1\over 2}\Delta(g)\]
where $\Delta$ extends canonically to $\Omega^1\tens_{C^\infty(M)}\Omega^1$. Although this kind of view of Ricci is known in specially adapted coordinates, this formula via the Hodge Laplacian $\Delta=\extd\delta+\delta\Delta$ appears to be new, see Section~2.5.

Next, in Section~3, we show how the above arises as a special case of a new algebraic point of view on metrics and connections as arising out of the purely algebraic problem of centrally extending a differential graded algebra by an additional graded-central 1-form $\theta'$ with $\extd\theta'=0$. We introduce the notion of a cleft extension here and show that cleft extensions of the classical exterior algebra $(\Omega(M),\extd)$ correspond to a class of metric-connection pairs where the metric-compatibiity tensor and torsion are matched (Proposition~\ref{cleftclass}). Within this theory of cleft extensions we consider those which are isomorphic to ones where $\extd$ is not changed, which we call `flat'. In the classical case this lands us on the classical Levi-Cevita connection from our new point of view. In the process we put the latter into a more general context where we now think of a metric-connection pair as equivalent to extension or `cocycle' data $(\Delta,\lo\ ,\ \lc)$ for the classical exterior algebra as a noncommutative differential graded algebra (DGA), i.e. as quantisation data for the differential structure. 

Much more than this, our algebraic approach works equally well with the initial DGA (the one being centrally extended) already being a non-graded-commutative or `quantum' $\Omega(A)$. This amounts then to a new approach to noncommutative Riemannian geometry as cleft extension data for $\Omega(A)$ (see Proposition~\ref{cleftcon1}) in the form of an associated quantum metric inner product $(\ ,\ ):\Omega^1(A)\tens_A\Omega^1(A)\to A$ and a quantum bimodule covariant derivative (as defined in the Preliminaries). Similar `bimodule connections' are by now much studied  in noncommutative geometry\cite{DV2,BegMa2,BegMa3} and it is remarkable that essentially such objects, albeit slightly generalised to covariant derivative along forms of all degrees and with $(\ ,\ )$ not necessarily non-degenerate, play the same as controlling the extensions of $\Omega(A)$. In the process we solve a couple of open problems in noncommutative geometry at least for the class of bimodule covariant derivatives arising here, specifically the notion of noncommutative or `quantum' inner product and in the invertible case a proposal for {\rm Ricci} that is well-defined when the quantum metric is quantum-symmetric in the sense of killed under the wedge product.  A main result at this level of general $\Omega(A)$ is Theorem~\ref{construct} which from a map $\perp$ obeying a 4-term identity
\[ (-1)^{|\eta|} (\omega\eta)\perp\zeta+(\omega\perp\eta)\zeta=\omega\perp(\eta\zeta)+(-1)^{|\omega|+|\eta|}\omega(\eta\perp\zeta),\ \forall\omega,\eta,\zeta\in\Omega\hfill\quad(**)\]
and compatible $\delta$ constructs a flat cleft central extension of the type which, in the classical case (see Theorem~\ref{cleftclasslevi}) is the Levi-Civita one.  Indeed, in the classical case (**) is solved automatically by extending the metric as a biderivation (see above) but in the general case this depends on the relations of $\Omega(A)$ and we instead specify $\perp$ at all degrees. Apart from this complication, Theorem~\ref{construct} plays the role of the Koszul formula in noncommutative geometry in giving the connection from the `metric' data $\perp$. The actual quantum metric is closely related to the latter and is constructed at the same time. We include an illustrative non-graded-commutative Example~\ref{Z2} on a set of two points. 

The  `homologically trivial' case where the DGA is isomorphic to a tensor product is also of interest and has recently been used to study stochastic differentials \cite{Beggs} cf. earlier work \cite{MH}. However, it is our very different cleft case which corresponds to Riemannian geometry and which we study. At the end of the Section~3, Proposition~\ref{Omegafull}, we consider a further extension problem where we also allow $\extd\theta'$, and here again (*) emerges as the solution to the algebraic extension problem. We denote the two extended DGAs by $\tilde{\Omega}(A)$ and $\tilde{\tilde{\Omega}}(A)$ respectively. 

Section~4 turns to a specific application but starts with another algebraic construction, namely the introduction of `time' by a graded-semidirect product of a DGA by the quantum DGA in one variable $\Omega(t,\extd t)$ (Proposition~\ref{semicalc}). The main result is the construction of an action of $\Omega(t,\extd t)$ by derivations on $\tilde{\tilde{\Omega}}(A)$ in the case where $A$ is functions on a smooth Riemannian manifold $M$. This comes down to further results in classical Riemannian geometry but in the approach of Section~2 based on the codifferential. Specifically, we introduce the notion that a 1-form $\tau$ is `$\delta$-conformal' if 
\[  [\delta,\CL_\tau]\omega=\alpha\delta\omega + (|\omega|-\beta)\inter_{\extd\alpha}\omega,\quad\forall \omega\in \Omega^1(M),\]
for some function $\alpha$ and some constant $\beta$. Lemma~\ref{metrictau} shows that in the classical setting this is equivalent to more conventional notions of conformal Killing 1-forms \cite{Yano}. Among necessary classical results (Proposition~\ref{tauprop}),  we prove that if the above condition holds in degree 1 as stated then it holds for all degrees of $\omega\in \Omega$ as long as $(\ ,\ )$ is nondegenerate. We also prove a similar commutation rule for $[\Delta,\CL_\tau]$ in all degrees. The resulting DGA  $\tilde{\tilde{\Omega}}(M)\lcross\Omega(t,\extd t)$ is a noncommutative differential version or `quantisation' of $M\times\R$ and extends  \cite[Sec.~3]{Ma:bh} to all degrees of forms although not quite  in the full generality used there to quantise spacetimes such as the  Schwarzschild black-hole. The general case in \cite[Sec.~3]{Ma:bh} should also be achievable by the same methods as here but with a cocycle semidirect product.  

The present paper does not depend on noncommutative geometry other than the use of differential graded algebras common to most approaches, including that of Connes\cite{Con}. The latter reconstructs a manifold from its Dirac operator and then takes that as a point of view for the noncommutative case, which is certainly more sophisticated than reconstructing from the codifferential or divergence as we do.  Further noncommutative geometric background in our case comes from the `quantum groups' approach where central extensions arise in the semisimple case as the leading correction to the bicovariant differential graded algebra\cite{BegMa1, MaTao}, and from the theory of (bimodule) connections defined algebraically \cite{DV2,BegMa2}. Although our starting differential graded algebra will tend to be of classical type,  Proposition~\ref{Omegaext} and Proposition~\ref{semicalc} are general and more fully noncommutative examples of them will be taken up elsewhere following on from recent results for inner DGAs in \cite{Ma:gra}. 

\subsection{Preliminaries} We take an algebraic approach in which we work with an algebra $A$ over a field $k$ of characteristic not 2. Our main application is to manifolds and so the field could be taken to be $\R$ and the algebra could be taken to be smooth functions on a smooth manifold. We require enough differentiable structure so as to have an associative  `differential graded algebra' (DGA) of differential forms,  $\Omega(A)=\oplus_n\Omega^n$ where $\Omega^0=A$, equipped with a graded-derivation $\extd:\Omega^i\to \Omega^{i+1}$ with $\extd^2=0$. We say that a DGA is {\em standard} if $\Omega^1$ is spanned by elements of the form $a\extd b$ for $a,b\in A$ and $\Omega$ is generated by degrees 0,1 over $A$.  We are mainly interested in the case of {\em classical type} where  $\Omega(A)$ is graded-commutative, standard, and given by the tensor algebra over $A$ of $\Omega^1$ modulo relations of antisymmetry. This is intended to keep us close to the classical situation and ensures in particular that antisymmetric module maps descend to $\Omega(A)$. 

We will always work with differential forms, but once we have an `metric inner product', by which we mean a bimodule map $(\ ,\ ):\Omega^1\tens_A\Omega^1\to A$, we will have an associated `vector field'  $(\omega,\ ):\Omega^1\to A$ or $X_\omega=(\omega,\extd(\ )):A\to A$ for any $\omega\in \Omega^1$ and this will be relevant to the motivation behind some of the definitions in the paper.  The meaning of $(\ ,\ )$ non-degenerate is the obvious one but one way to achieve it is the existence of a central element $g\in\Omega^1\tens_A\Omega^1$ `the metric' such that $(\omega, g\uo)g\ut=\omega=g\uo(g\ut,\omega)$ for all $\omega\in\Omega^1$. Here $g=g\uo\tens g\ut$ (a sum of such terms understood) is a notation.  This is the normal set-up in noncommutative differential geometry but is also useful in the `classical' case, where we  normally also require that $(\ ,\ )$ is symmetric.

By a (left) algebraic connection on a DGA in noncommutative geometry one normally means $\nabla:\Omega^1\to \Omega^1\tens_A\Omega^1$ in degree 1 or more generally $\Omega^m\to \Omega^1\tens_A\Omega^m$ such that $\nabla(a\eta)=a\nabla\eta+\extd a\tens\nabla\eta$ for all $a\in A$ and $\eta\in \Omega^m$.  The nicest case is that of a `bimodule connection' where in addition we have $\nabla(\eta a)=(\nabla\eta)a+\sigma(\omega\tens_A\extd a)$ for some map $\sigma:\Omega^m\tens_A\Omega^1\to \Omega^1\tens_A\Omega^m$, called the `generalised braiding'. Such a map if it exists is uniquely determined, so this is a really a property that a left connection can further have. 

One departure, we shall more often be interested in directly defining a `1-form covariant derivative' $\nabla_\omega:\Omega^m\to \Omega^m$ for all $\omega\in \Omega^1$ with analogous properties given by evaluation against a map $(\ ,\ ):\Omega^1(A)\tens_A\Omega^1(A)\to A$, namely
\begin{equation}\label{leftcov} \nabla_{a\omega}=a\nabla_\omega,\quad \nabla_\omega(a\eta)=\nabla_{\omega a}\eta+ (\omega,\extd a)\eta,\quad\forall \omega\in \Omega^1,\ \eta\in\Omega^m .\end{equation}
and this is a bimodule covariant derivative if 
\begin{equation}\label{leftbimod} \nabla_\omega(\eta a)=(\nabla_\omega\eta)a+\sigma_\omega(\eta\tens_A\extd a);\quad \sigma:\Omega^1\tens_A\Omega^m\tens_A\Omega^1\to \Omega^m\end{equation} 
for some bimodule map $\sigma$. Again this is a property of a covariant derivative rather than additional structure. Also for any covariant derivative we have the `half curvature'
\[  \rho(\omega\tens_A\eta)=\nabla_\omega\nabla_\eta-\nabla_{\nabla_\omega\eta},\quad \forall \omega,\eta\in \Omega^1\]
and it is a nice check from the above properties that this depends only on $\omega\tens_A\eta$: 
\begin{eqnarray*}
\nabla_{\omega}\nabla_{a\eta}-\nabla_{\nabla_{\omega}(a\eta)}&=& \nabla_{\omega}(a\nabla_\eta\ )-\nabla_{\nabla_{\omega a}\eta+(\omega,\extd a)\eta}\\
&=&\nabla_{\omega a}\nabla_{\eta}-\nabla_{\nabla_{\omega a}\eta}+(\omega,\extd a)\nabla_\eta-\nabla_{(\omega,\extd a)\eta}\\
&=&\nabla_{\omega a}\nabla_{\eta}-\nabla_{\nabla_{\omega a}\eta}
\end{eqnarray*}
As a result when $(\ ,\ )$ is invertible the `Laplace-Beltrami' operator
\begin{equation}\label{LB}\Delta_{LB}=\nabla_{g\uo}\nabla_{g\ut}-\nabla_{\nabla_{g\uo}g\ut}\end{equation}
is well-defined. It is not intended to go deeply into noncommutative geometry but some of our constructions will be no harder in the possibly noncommutative case and this is one of them.

We will often be interested in the failure of the Leibniz rule. For this it is usual to define for any degree $d$ linear map $B:\Omega(A)\to \Omega(A)$, the `Leibnizator' 
\[ L_B(\omega,\eta)=B(\omega\eta)-(B\omega)\eta-(-1)^{d|\omega|}\omega B\eta,\quad \forall \omega,\eta\in\Omega(A).\]

\section{Reconstruction from a codifferential}

In the case of a Riemannian manifold (or pseudo-Riemannian, of any signature) there is a Hodge $\star$-operator and from this a codifferential $\delta:\Omega^i\to \Omega^{i-1}$. It will be immediately clear that we can recover the Riemannian structure from $\delta$ because we can recover the metric inner product $(\ ,\ )$ according to Definition~\ref{regular}. However, this point of view turns out to give natural formulae for all of the ensuing structures and these formulae are exactly what are needed for the quantisation in Section~3.

\subsection{Interior product}

\begin{definition} \label{regular} Let $\Omega(A)$ be a standard DGA. We say that a degree -1 linear map is {\em regular} if there exist degree -1 bimodule maps $\inter:\Omega^1\tens_A\Omega\to \Omega$ and $\rinter:\Omega\tens_A\Omega^1\to \Omega$ such that 
\[ \delta(a\omega)=a\delta\omega+\inter_{\extd a}\omega,\quad \delta(\omega a)=(\delta\omega)a+\omega\rinter_{\extd a}\]
where $\rinter$ acts from the right.  In this case we refer to the associated bimodule map $(\ ,\ ):\Omega^1\tens_A\Omega^1\to A$ defined by
\[ (\omega,\eta)={1\over 2}(\omega\rinter_\eta+\inter_\omega\eta),\quad \omega,\eta\in \Omega^1\]
as the associated `metric inner product'. \end{definition}
Note that if these maps exist, they are uniquely determined by $\delta$. 
\begin{lemma}\label{deriv} Let $\delta$ be regular. 

(1) $\delta$ anticommutes with $\inter_\eta,\rinter_\eta$ for all $\eta\in \Omega^1$ iff $\delta^2$ is a bimodule map. In this case $\inter_\eta,\rinter_\omega$ mutually anticommute for all $\eta,\omega\in \Omega^1$.

(2)  $\inter_{\extd a},\rinter_{\extd a}$ are graded-derivations iff
\[ L_\delta(a\omega,\eta)=a L_\delta(\omega,\eta)+(-1)^{|\omega|}\omega\inter_{\extd a}\eta,\quad L_\delta(\omega,\eta a)=L_\delta(\omega,\eta)a+(\omega\rinter_{\extd a})\eta,\quad\forall\omega,\eta\in \Omega,\ a\in A.\]
\end{lemma}
\proof (1) For any $a\in A, \omega\in\Omega$,  $\delta\inter_{\extd a}\omega=\delta(\delta(a\omega)-a\delta\omega)=\delta^2(a\omega)-a\delta^2\omega-\inter_{\extd a}\delta\omega$ and similarly for $\rinter$. When quasi-nilpotency holds, we then have
\begin{eqnarray*}\inter_{\extd a}(\omega\rinter_{\extd b})&=&\inter_{\extd a}(\delta(\omega b)-(\delta\omega)b)=-\delta((\inter_{\extd a}\omega)b)+(\delta\inter_{\extd a}\omega)b=-(\inter_{\extd a}\omega)\rinter_{\extd b}.\end{eqnarray*}
This implies that $(\inter_{(\extd a)b}\omega)\rinter_{f\extd h}=(\inter_{\extd a}b\omega f)\rinter_{\extd h}=-\inter_{\extd a}(b\omega f\rinter_{\extd h})=-\inter_{(\extd a)b}(\omega \rinter_{f\extd h})$ for all $a,b,f,h\in A$ by the bimodule map properties, hence the result applies to arbitrary 1-forms.  (2) We also have
\begin{eqnarray*}L_\delta(a\omega,\eta)-aL_\delta(\omega,\eta)=\delta(a\omega\eta)-(\delta(a\omega))\eta-a\delta(\omega\eta)+a(\delta\omega)\eta=\inter_{\extd a}(\omega\eta)-(\inter_{\extd a}\omega)\eta\end{eqnarray*}
for all $\omega,\eta\in \Omega$. Hence $\inter_{\extd a}$ is a right derivation iff the first stated condition holds. Similarly for $\rinter_{\extd a}$ (as a left derivation). The conditions amount to 6-term conditions on the behaviour of $\delta$ on a triple product where one factor is in $A$. 
\endproof

In the graded-commutative case the left and right `interior products' coincide. It follows if we assume that $\delta^2$ is a tensorial (a module map) that $\inter_\omega^2=0$ for all $\omega\in\Omega^1(A)$. Similarly, the derivation property extends to $\inter_\omega$ for all $\omega\in\Omega^1(A)$. It is also convenient for contact with classical differential geometry (but not essential as we will see in Section~3) to suppose that 
\begin{equation}\label{symminter}\inter_{\extd a}(\extd b)=\inter_{\extd b}(\extd a),\quad\forall a,b\in A.\end{equation}
 Accordingly:

\begin{definition}\label{classicaldelta} We say that $\delta$ is of {\em classical type} if it is regular, the two conditions of Lemma~\ref{deriv} apply and symmetry in the form (\ref{symminter}) holds. When both $\Omega(A)$ and $\delta$ are of classical type we say that the pair $(\Omega(A),\delta)$ is of classical type.
\end{definition}

 In this case the anticommutativty of $\inter$ means that we can extend it to $\inter_{\omega_1\cdot \omega_m}=\inter_{\omega_1}\cdots\inter_{\omega_m}$ where $\omega_i\in\Omega^1$, to give a well-defined degree $-m$ linear map on $\Omega(A)$. Note, however, that this map is in general only a graded-derivation when $m=1$. For example one may easily compute 
\begin{equation}\label{i2 lie} L_{\inter_{\omega_1 \omega_2}}(\omega,\eta)=(-1)^{|\omega|}\left((\inter_{\omega_1}\omega)\inter_{\omega_2}\eta-(\inter_{\omega_2}\omega)\inter_{\omega_1}\eta\right),\quad\forall\omega,\eta\in\Omega\end{equation} and similar formulae in general. Using this notation, the mutual anticommutativity of $\inter_{\extd a},\delta$ in the proof of Lemma~\ref{deriv} is readily seen to generalise to
 \begin{equation}\label{delta lie} \delta\inter_\omega + \inter_\omega\delta=\inter_{\extd\omega},\quad\forall \omega\in\Omega^1.\end{equation}
(This implies a similar formulae for all degrees of $\omega$ but with a graded-commutator on the left and side.)

\begin{lemma}\label{Schouten} Let $(\Omega(A),\delta)$ be of classical type. Then
\[\inter_\zeta L_\delta(\omega,\eta)=-L_\delta(\inter_\zeta\omega,\eta)-(-1)^{|\omega|}L_\delta(\omega,\inter_\zeta\eta)+L_{\inter_{\extd\zeta}}(\omega,\eta),\quad \forall\omega,\eta\in \Omega,\ \zeta\in \Omega^1\]
and \[ \inter_\zeta L_\delta(\omega,\eta)=\inter_\omega \extd\inter_\eta\zeta-\inter_\eta\extd\inter_\omega\zeta-\inter_\eta\inter_\omega\extd\zeta,\quad\forall\omega,\eta,\zeta\in\Omega^1.\]
If, moreover,  $(\ ,\ )$ is nondegenerate then $L_\delta(\omega,\ )$ is a degree $|\omega|-1$ derivation for all $\omega\in \Omega$ and 
\[  L_\delta(\omega_1\cdots\omega_m,\eta_1\cdots\eta_n)=\sum_{i,j}(-1)^{i+j}\omega_1\cdots\widehat{\omega_i}\cdots\omega_m L_\delta(\omega_i,\eta_j)\eta_1\cdots\widehat{\eta_j}\cdots\eta_n\]
where $\omega_i,\eta_j\in \Omega^1$. Here the hat denotes omission.
\end{lemma}
\proof (1) We use (\ref{delta lie}) in the definition of $L_\delta$. Thus
\begin{eqnarray*} \inter_\zeta L_\delta(\omega,\eta)&=&\inter_\zeta\delta(\omega\eta)-(\inter_\zeta\delta\omega)\eta-(-1)^{|\omega|-1}(\delta\omega)\inter_\zeta\eta-(-1)^{|\omega|}(\inter_\zeta\omega)\delta\eta-\omega\inter_\zeta\delta\eta\\
&=&-\delta((\inter_\zeta\omega)\eta+(-1)^{|\omega|}\omega\inter_\zeta\eta)+\inter_{\extd\zeta}(\omega\eta)+(\delta\inter_\zeta\omega)\eta-(\inter_{\extd\zeta}\omega)\eta\\
&&\quad+(-1)^{|\omega|}(\delta\omega)\inter_\zeta\eta-(-1)^{|\omega|}(\inter_\zeta\omega)\delta\eta+\omega\delta\inter_\zeta\eta -\omega\inter_{\extd\zeta}\eta\\
&&=-L_\delta((\inter_\zeta\omega),\eta)-(-1)^{|\omega|}L_\delta(\omega,\inter_\zeta\eta)+\inter_{\extd\zeta}(\omega\eta)-(\inter_{\extd\zeta}\omega)\eta -\omega\inter_{\extd\zeta}\eta
\end{eqnarray*}
We also have $L_\delta(a,\omega)=L_\delta(\omega,a)=\inter_{\extd a}\omega$ for all $a\in A$ and $\omega\in\Omega$ after which our result implies the explicit formula in the case $\omega,\eta\in\Omega^1$. (1) Note that $L_{\inter_{\extd\zeta}}(\omega, )$ using (\ref{i2 lie}) is a graded-derivation of degree $\omega$ and when $\omega\in \Omega^1$ our stated formula implies 
\begin{equation}\label{lie omega1} \inter_\zeta L_\delta(\omega,\eta)=-\inter_{\extd\inter_\zeta\omega}(\eta)+L_\delta(\omega,\inter_\zeta(\eta))+L_{\inter_{\extd\zeta}}(\omega,\eta).\end{equation}
We prove by induction that $L_\delta(\omega,\ )$ is a derivation on a product $\eta\eta'$, assuming that the same is true on a product where either $\eta,\eta'$ are replaced by a form of one less degree. Thus
\begin{eqnarray*}\inter_\zeta L_\delta(\omega,\eta\eta')&=&-\inter_{\extd\inter_\zeta\omega}(\eta\eta')+L_\delta(\omega,\inter_\zeta(\eta\eta'))+L_{\inter_{\extd\zeta}}(\omega,\eta\eta')\\
&=&-(\inter_{\extd\inter_\zeta\omega}\eta)\eta'-(-1)^{|\eta|}\eta\inter_{\extd\inter_\zeta\omega}(\eta')+L_{\inter_{\extd\zeta}}(\omega,\eta)\eta'+(-1)^{|\eta|}\eta L_{\inter_{\extd\zeta}}(\omega,\eta)\\
&&+L_\delta(\omega,(\inter_\zeta\eta))\eta'+(\inter_\zeta\eta)L_\delta(\omega,\eta')+(-1)^{|\eta|}L_\delta(\omega,\eta)\inter_\zeta\eta'
+(-1)^{|\eta|}\eta L_\delta(\omega,\inter_\zeta\eta')\\
&=&(\inter_\zeta\eta)L_\delta(\omega,\eta')+(-1)^{|\eta|}\eta\inter_{\zeta}L_\delta(\omega,\eta')+(\inter_\zeta L_\delta(\omega,\eta))\eta'+(-1)^{|\eta|}L_\delta(\omega,\eta)\inter_\zeta\eta'\\
&=&\inter_\zeta\left(L_\delta(\omega,\eta)\eta'+\eta L_\delta(\omega,\eta')\right)
\end{eqnarray*}
using (\ref{lie omega1}) on the product and in reverse to recognise the answer. Now, if $\eta$ or $\eta'$ have degree 0 then the derivation property on $\eta\eta'$ reduces to part of Lemma~\ref{deriv}, so this holds and provides the boundary condition for the induction. Thus if $(\ ,\ )$ is nondegenerate we see that $L_\delta(\omega,\ )$ is a derivation for all $\omega\in \Omega^1$. We now prove that if the DGA is graded-commutative and $L_\delta(\omega, )$ is a derivation for $\omega\in\Omega^1$  then $L_\delta(\omega,\ )$ is a graded-derivation of degree $|\omega|-1$ for all $\omega$. The degree zero case $L_\delta(a,\ )=\inter_{\extd a}$ is already assumed to be a graded-derivation of degree -1 in Lemma~\ref{deriv}. We will need the tautological identity
\begin{equation}\label{L1} L_\delta(\omega\eta,\zeta)+L_\delta(\omega,\eta)\zeta=L_\delta(\omega,\eta\zeta)+(-1)^{|\omega|}\omega L_\delta(\eta,\zeta)\end{equation}
which holds for the Leibnizator of any degree -1 linear map on any graded algebra (just write out the definitions on both sides), cf\cite{Coll} in the graded-commutative case. Suppose  $L_\delta(\omega,\ )$ is a degree $|\omega|-1$ derivation for $\omega$ of some degree. Using (\ref{L1}) we deduce
\[ L_\delta(\omega\eta,\zeta)=L_\delta(\omega,\eta\zeta)+(-1)^{|\omega|}\omega L_\delta(\eta,\zeta)-L_\delta(\omega,\eta)\zeta\]
\begin{equation}\label{L3}\quad\quad\quad=(-1)^{(|\omega|-1)|\eta|}\eta L_\delta(\omega,\zeta)+(-1)^{|\omega|}\omega L_\delta(\eta,\zeta)\end{equation}
for the given $\omega$ and all $\eta,\zeta$. We also suppose $L_\delta(\eta,\ )$ is a degree $|\eta|-1$ graded derivation for $|\eta|\le |\omega$ (it suffices to take $|\eta|=1$). Then,
\begin{eqnarray*} 
L_\delta(\omega\eta,\xi\zeta)&=&L_\delta(\omega,\eta\xi\zeta)+(-1)^{|\omega|}\omega L_\delta(\eta,\xi\zeta)-L_\delta(\omega,\eta)\xi\zeta\\
&=&(-1)^{(|\omega|-1)(|\eta|+|\xi|)}\eta\xi L_\delta(\omega,\zeta)+L_\delta(\omega,\eta\xi)\zeta+(-1)^{|\omega|}(-1)^{(|\eta|-1)|\xi|}\omega \xi L_\delta(\eta,\zeta)\\
&&+(-1)^{|\omega|}\omega L_\delta(\eta,\xi)\zeta-L_\delta(\omega,\eta)\xi\zeta\\
&=&(-1)^{(|\omega|-1)(|\eta|+|\xi|)}\eta\xi L_\delta(\omega,\zeta)+(-1)^{|\omega|}(-1)^{(|\eta|-1)|\xi|}\omega \xi L_\delta(\eta,\zeta)+L_\delta(\omega\eta,\xi)\zeta\\
&=&(-1)^{(|\omega|+|\eta|-1)|\xi|}\left( (-1)^{(|\omega|-1)|\eta|}\xi\eta L_\delta(\omega,\zeta)+(-1)^{|\omega|} \xi \omega L_\delta(\eta,\zeta)\right)+L_\delta(\omega\eta,\xi)\zeta\\
&=&(-1)^{(|\omega|+|\eta|-1)|\xi|}\xi L_\delta(\omega\eta,\zeta)+L_\delta(\omega\eta,\xi)\zeta
\end{eqnarray*}
where we used (\ref{L1}), then our assumed derivation properties, (\ref{L1}) in reverse, graded-commutativity and the computation above, to recognise the answer. This proves the required graded-derivation property by induction. It follows from (\ref{L3}) in the graded-commutative case that if a degree -1 bilinear map $L_\delta$ obeys (\ref{L1}) and $L_\delta(\omega,\ )$ is a graded-derivation, then
\begin{equation}\label{Lmulti} L_\delta(\omega_1\cdots\omega_m,\ )=\sum_{i=1}^m(-1)^{i-1}\omega_1\cdots\widehat{\omega_i}\cdots\omega_m L_\delta(\omega_i,\ ),\quad\forall\omega_i\in \Omega^1\end{equation}
leading to the formula stated. This is a general observation which we will also use for other maps obeying  (\ref{L1}). 
 \endproof

The specific formula for $\omega,\eta$ of degree 1 confirms that $L_\delta(\omega,\eta)$  in the case of a classical manifold corresponds via the metric to the Lie bracket of vector fields. Thus, let  $X_\omega=(\omega,\extd(\ ))$ be the vector field associated to a 1-form and let $[X_\omega,X_\eta]$  be the usual Lie bracket of such vector fields viewed as a tensorial map on 1-forms. Then
\begin{equation}\label{liebra} [X_\omega,X_\eta](\zeta)=(\omega,\extd(\eta,\zeta))-(\extd(\omega,\zeta),\eta)-\inter_\eta\inter_\omega\extd\zeta,\quad\forall \omega,\eta,\zeta\in \Omega^1\end{equation}
in agreement with $\inter_\zeta L_\delta(\omega,\eta)$ in the lemma. In this case it is clear from the formula on higher degrees that  $L_\delta(\omega,\eta)$ corresponds to the Schouten bracket of alternating multivector fields cf. \cite{Coll}, and indeed the results of  Lemma~\ref{Schouten} can be seen as parallel to properties of this \cite{Sch,Marle}. There will also be a form of graded Jacobi identity which we have not elaborated here as we will not need it. On the other hand we view $L_\delta(\omega,\eta)$ as the primary object with its evaluations such as (\ref{liebra}) defining the bracket as a linear map on $\zeta$ of appropriate degree even in the degenerate case. For example, by iterating Lemma~\ref{Schouten} one has
\begin{equation}\label{liehigher} \inter_{\zeta}L_\delta(\omega,\eta)=\inter_\omega\extd\inter_\eta \zeta-\inter_\eta\extd\inter_\omega\zeta-\inter_{\eta}\inter_\omega\extd\zeta,\quad\forall \omega\in \Omega^1,\ \eta,\zeta\in\Omega^2\end{equation}
without assuming nondegeneracy, and similarly in general degree for $\eta,\zeta$ (we will need this only in degrees 1,2). 

\begin{corollary}\label{delta triple} If $(\Omega(A),\delta)$ is of classical type and $(\ ,\ )$ is nondegenerate then
\[ \delta(\omega\eta\zeta)=(\delta(\omega\eta))\zeta+(-1)^{|\omega|}\omega\delta(\eta\zeta)+(-1)^{(|\omega|-1)|\eta|}\eta\delta(\omega\zeta)\]
\[\qquad\qquad- (\delta\omega)\eta\zeta-(-1)^{|\omega|}\omega(\delta\eta)\zeta-(-1)^{|\omega|+|\eta|}\omega\eta\delta\zeta\]
$\forall \omega,\eta,\zeta\in \Omega$. In other words, $(\Omega(A),\delta)$ is a Batalin-Vilkovisky algebra slightly generalised to allow $\delta^2$ to be a left module map. \end{corollary}
\proof This is just the content of $L_\delta(\omega,\ )$ a graded-derivation, written out in terms of $\delta$.    
 \endproof
The case of $\omega$ of degree 0 is the content of our classical type assumption (the derivation properties in Lemma~\ref{deriv}) and our result is that in the nondegenerate case the stated identity then holds in all degrees. 

Also in the case of classical type we now introduce an operation $\perp: \Omega\tens_A\Omega\to \Omega$ of degree -2,
\[ (\omega_1\cdots \omega_m)\perp (\eta_1\cdots\eta_n)=\sum_{i,j}(-1)^{i+j}(\omega_i,\eta_j)\omega_1\cdots\widehat{\omega_i}\cdots \omega_m \eta_1\cdots \widehat{\eta_j}\cdots \eta_n,\quad\omega_i,\eta_j\in \Omega^1.\]
Classically $\inter_\omega$ in degree 1 is a graded derivation of degree -1 and in our case similarly
\[  \omega_1\cdots\omega_m\perp(\ )=\sum_{j=1}^m (-1)^{j-1} \omega_1\cdots\widehat{\omega_j}\cdots\omega_m \inter_{\omega_j}(\ ),\]
is a degree $m-2$ derivation, while
\[  (\ )\perp\eta_1\cdots\eta_n=\sum_{i=1}^n (-1)^{i-1} (\inter_{\eta_i}\ )\eta_1\cdots\widehat{\eta_i}\cdots\eta_n\]
is such that $((-1)^D(\ ))\perp\eta_1\cdots\eta_n$ is a degree $n-2$ right derivation, where $D$ is the degree operator. It is also clear from the form of these expressions that they depend tensorialy (they are $A$-module maps) and antisymmetrically and hence descend to $\Omega(A)$. In particular, if $\omega$ is degree 1 then $\omega\perp$ and $\perp\omega$ revert to the interior product by $\omega$. We will particularly need 
\[ \omega_1\omega_2\perp= \omega_2\inter_{\omega_1}- \omega_1\inter_{\omega_2}.\]
Note that interior products are usually considered by vector fields and one could consider that the 1-form $\omega$ is being converted to a vector field  $(\omega,\extd(\ ))$ for this purpose. Later on, in Section~3, we shall generalise this construction so that $\perp$ need not be symmetric when restricted to degree 1, but for $\delta$ of classical type as here, $\perp$ just extends the metric inner product $(\ ,\ )$. 

Using the interior product we define the  `Lie derivative'  as 
\[ \CL_\omega=\extd\inter_\omega+\inter_\omega\extd,\quad\omega\in\Omega^1\]
along the lines of the classical Cartan formula. Clearly it obeys
\begin{equation}\label{lie mod} \CL_{a\omega}\eta=a\CL_\omega\eta+(\extd  a)\inter_\omega\eta,\quad \CL_\omega(a\eta)=a\CL_\omega(\eta)+ (\omega,\extd a)\eta,\quad\forall\omega\in\Omega^1.\end{equation}

\subsection{Form covariant derivatives}

We similarly define a `covariant derivative'  $\nabla:\Omega^1\times\Omega\to \Omega$ to be a map characterised by (\ref{leftcov}) as explained in the Preliminaries. We will be interested in metric compatibility, which means vanishing of the tensor 
\begin{equation}\label{metriccomp} C_\omega(\eta,\zeta)=(\omega,\extd(\eta,\zeta))-(\nabla_\omega\eta,\zeta)-(\eta,\nabla_\omega\zeta),\quad\forall \omega,\eta,\zeta\in \Omega^1.\end{equation}
In the same vein we define
\begin{equation}\label{T} T(\omega,\eta)(\zeta)=(\omega,\nabla_\eta\zeta)-(\eta,\nabla_\omega\zeta)-\inter_{\omega}\inter_\eta\extd\zeta,\quad\forall \omega,\eta,\zeta\in \Omega^1\end{equation}
as the torsion of a covariant derivative. Both maps are easily seen to be tensorial in all of their inputs. These formulae are dualisations of the usual formulae with vector fields and make sense in this form for any standard graded commutative $\Omega(A)$ of classical type, but note that we do not assume that $(\ ,\ )$ is nondegenerate.

\begin{lemma}\label{symanti} Let $\Omega(A)$ be of classical type and equipped with a symmetric metric inner product $(\ ,\ )$ and  $\nabla$ a covariant derivative. Then
\[T(\omega,\eta)(\zeta)+C_\eta(\omega,\zeta)-C_\omega(\eta,\zeta)=\inter_\zeta\left(\nabla_\omega\eta-\nabla_\omega\eta\right)-\inter_\omega \extd\inter_\eta\zeta+\inter_{\eta}\extd\inter_\omega\zeta-\inter_\omega\inter_\eta\extd\zeta\]
\[T(\zeta,\omega)(\eta)+T(\zeta,\eta)(\omega)-C_\zeta(\omega,\eta)=  \inter_\zeta\left(\nabla_\omega\eta+\nabla_\eta\omega-\CL_\omega\eta-\CL_\eta\omega+\extd(\omega,\eta)\right)\]
for all $\omega,\eta,\zeta\in \Omega^1$.
\end{lemma}
\proof For the first part we use (\ref{metriccomp}) in each of the first two terms of the definition (\ref{T}) of torsion. For the second part we use (\ref{T}) on each term to compute
\begin{eqnarray*} (\zeta,\nabla_\omega\eta&+&\nabla_\eta\omega)\\
&&=(\omega,\nabla_\zeta\eta)+T(\zeta,\omega)(\eta)+ \inter_\zeta\inter_\omega\extd\eta+ (\eta,\nabla_\zeta\omega)+T(\zeta,\eta)(\omega)+ \inter_\zeta\inter_\eta\extd\omega\\
&&=T(\zeta,\omega)(\eta)+T(\zeta,\eta)(\omega)+\inter_\zeta(\extd(\omega,\eta)+\inter_\omega\extd\eta+\inter_\eta\extd\omega)-C_\zeta(\omega,\eta)
\end{eqnarray*}
using (\ref{metriccomp}). We then recognise the answer in terms of a Lie derivative. 
\endproof

This means that in the nondegenerate case a metric compatible torsion free covariant derivative, if it exists, is uniquely determined as its symmetric and antisymmetric parts are determined (as on a classical manifold). One can also treat the curvature evaluated against one-forms in a similar spirit. However, in the case where $(\ ,\ )$ comes from a $\delta$ of classical type one can do rather better:

\begin{proposition}\label{curvT} Let $(\Omega(A),\delta)$ be of classical type and $\nabla$ a covariant derivative. Then
\[ R(\omega,\eta)(\zeta):=\nabla_\omega\nabla_\eta\zeta- \nabla_\eta\nabla_\omega\zeta-\nabla_{L_\delta(\omega,\eta)}\zeta,\quad \forall \omega,\eta,\zeta\in\Omega^1\]
is tensorial in all its inputs and reduces to the usual curvature in the nondegenerate or algebraic cases. If the covariant derivative is $(\ ,\ )$-compatible then
\[ T(\omega,\eta):=\nabla_\omega\eta-\nabla_\eta\omega-L_\delta(\omega,\eta), \quad \forall \omega,\eta\in\Omega^1\]
is tensorial in its inputs and evaluates via $(\ ,\ )$ to the torsion. 
\end{proposition}
\proof We now let $\nabla_\omega$ be any covariant derivative and check
\begin{eqnarray*} R(\omega,\eta)(a\zeta)&=&\nabla_\omega(a\nabla_\eta\zeta)-\nabla_\eta(a\nabla_\omega\zeta)-a\nabla_{L_\delta(\omega,\eta)}\zeta\\
&&+\nabla_\omega((\eta,\extd a)\zeta)-\nabla_\eta((\omega,\extd a)\zeta)-(L_\delta(\omega,\eta),\extd a)\zeta\\
&=&a R(\omega,\eta)(\zeta)+(\omega,\extd(\eta,\extd a))-(\eta,\extd(\omega,\extd a))-(L_\delta(\omega,\eta),\extd a)=aR(\omega,\eta)(\zeta)\\
R(a\omega,\eta)(\zeta)&=&a\nabla_\omega\nabla_\eta\zeta-\nabla_\eta(a\nabla_\omega\zeta)-a\nabla_{L_\delta(\omega,\eta)}\zeta-(\eta,\extd a)\nabla_\omega\zeta=aR(\omega,\eta)(\zeta)
\end{eqnarray*}
for all $\omega,\eta,\zeta\in \Omega^1$, $a\in A$. For the first computation we used the defining property (\ref{leftcov}) of a covariant derivative followed by Lemma~\ref{Schouten}. For the second computation we used the covariant derivative property and Lemma~\ref{deriv}.  Similarly for the other input of the curvature. Note that at least in the nondegenerate case one can then evaluate the algebraic curvature $R_\nabla=(\extd\tens\id - (\wedge\tens\id)(\id\tens\nabla))\nabla$ by applying $\inter_\eta\inter_\omega$, to obtain
\[ R(\omega,\eta)(\zeta)=\nabla_\omega\nabla_\eta\zeta- \nabla_\eta\nabla_\omega\zeta-\nabla_{[\omega,\eta]}\zeta,\quad \forall \omega,\eta,\zeta\in\Omega^1\]
as a definition in this case, where the `Lie bracket' on 1-forms is given by  $[\omega,\eta]=L_\delta(\omega,\eta)$ or rather by its evaluation on 1-forms as explained above.  That we have $\inter_\eta\inter_\omega R_\nabla$ in the case where $\nabla:\Omega^1\to \Omega^1\tens_A\Omega^1$ is defined is part of the standard derivation of the algebraic expression for $R_\nabla$. Suffice for completeness to note by the covariant derivative property (\ref{leftcov}) that
\[ \nabla(\nabla_\omega\zeta)=\extd(\omega,\nabla^1\zeta)\tens\nabla^2\zeta+(\omega,\nabla^1\zeta)\nabla\nabla^2\zeta\]
where we use a notation $\nabla\zeta=\nabla^1\zeta\tens\nabla^2\zeta$. We use this in the evaluation of the 2nd term of $R_\nabla$ and Lemma~\ref{Schouten} for the evaluation of the first term $\inter_\eta\inter_\omega\extd$. 

The torsion $\nabla$ on 1-forms is likewise given by $\inter_\eta\inter_\omega$ against the algebraic torsion $T=\wedge\nabla-\extd:\Omega^1\to \Omega^2$ to give the map (\ref{T}). In the important case where the covariant derivative is $(\ ,\ )$-compatible, we recognise the formula in Lemma~\ref{Schouten} for $\inter_\zeta L_\delta(\omega,\eta)$ in the expression for torsion in Lemma~\ref{symanti}. We then take $T(\omega,\eta)\in \Omega^1$ as a definition applicable in the $(\ ,\ )$-compatible case.  Tensoriality is from Lemma~\ref{deriv} and (\ref{leftcov}). \endproof

\subsection{Levi-Civita covariant derivative}

We are now ready to state and prove our main result:

\begin{theorem}\label{levi} In the setting above with $(\Omega(A),\delta)$ of classical type, there is a covariant derivative 
\[ \nabla_{\omega}\eta={1\over 2}\left(L_\delta(\omega,\eta)+\CL_{\omega}\eta+(\extd\omega)\perp \eta\right),\quad\forall \omega\in\Omega^1,\ \eta\in\Omega.\]
which is torsion free and compatible with $(\ ,\ )$ and in the nondegenerate case acts as a derivation. 
\end{theorem}
\proof We will see in Section~3 how this formula arises as a requirement for quantisation; for the moment we verify directly that the stated expression is indeed a covariant derivative with the stated properties. Thus
\begin{eqnarray*} 2\nabla_{a\omega}\eta&=&\CL_{a\omega}(\eta)+L_\delta(a\omega,\eta)+(\extd (a\omega))\perp\eta\\
&=&a\CL_\omega\eta+(\extd a)\inter_{\omega}\eta+aL_\delta(\omega,\eta)-\omega\inter_{\extd a}\eta+((\extd a)\omega)\perp\eta+a\extd\omega\perp\eta\\
&=& a 2\nabla_\omega\eta\\
2\nabla_\omega(a\eta)&=&\CL_\omega(a\eta)+L_\delta(\omega,a\eta)+(\extd \omega)\perp a\eta\\
&=& a\CL_\omega\eta+(\omega,\extd a)\eta+aL_\delta(\omega,\eta)+(\inter_{\extd a}\omega)\eta +a(\extd\omega)\perp\eta=2\nabla_\omega(\eta)+2(\omega,\extd a).
\end{eqnarray*}

Next we note that acting in degree 1 this covariant derivative can be written as
\[ \nabla_\omega\eta={1\over 2}L_\delta(\omega,\eta)+ {1\over 2}(\CL_\omega\eta+\CL_\eta\omega - \extd(\omega,\eta))\]
so that
\begin{equation}\label{connpm} \nabla_\omega\eta-\nabla_\eta\omega=L_\delta(\omega,\eta),\quad \nabla_{\omega}\eta+\nabla_\eta\omega=\CL_{\omega}\eta+\CL_\eta\omega-\extd(\omega,\eta),\quad\forall \omega,\eta\in \Omega^1\end{equation}
which comparing with Lemma~\ref{symanti} and using the formula for $\inter_\zeta L_\delta(\omega,\eta)$ in Lemma~\ref{Schouten} implies that if $T=0$ then $C=0$.

It remains to prove that $T=0$ in (\ref{T}). For this we put in the particular form of our covariant derivative as found above. Then
\begin{eqnarray*} 2T(\omega,\eta)(\zeta)&=&(\omega,L_\delta(\eta,\zeta)-\CL_\eta\zeta+\inter_\zeta\extd\eta+2\extd\inter_\zeta\eta) - (\eta, \CL_\omega\zeta+L_\delta(\omega,\zeta)+\inter_\zeta\extd\omega)\\
&=&(\omega,L_\delta(\eta,\zeta))-(\eta,L_\delta(\omega,\zeta))+\inter_\omega\inter_\zeta\extd\eta+\inter_\omega\extd\inter_\zeta\eta-\inter_\eta\extd\inter_\omega\zeta-\inter_\eta\inter_\zeta\extd\omega\\
&=&-\inter_{\extd(\omega,\eta)}\eta+\inter_{\extd\omega}(\eta\zeta)+\inter_{\extd(\eta,\zeta)}\omega-\inter_{\extd\eta}(\omega\zeta)+\inter_\omega\inter_\zeta\extd\eta+\inter_\omega\extd\inter_\zeta\eta-\inter_\eta\extd\inter_\omega\zeta\\
&&\quad\quad-\inter_\eta\inter_\zeta\extd\omega\\
&=&0.
\end{eqnarray*}
We used Lemma~\ref{Schouten} for $L_\delta$ and symmetry of same-degree interior products to cancel.  

The Lie derivative and $\extd\omega\perp$ act by derivations which covers the symmetric part, and in the nondegenerate case  Lemma~\ref{Schouten} tells us that the antisymmetric part also acts as a derivation.  \endproof

This provides a natural `Levi-Civita' covariant derivative in our setting. In principle its curvature and other geometrical properties can be computed in terms of $\delta$. We also at the same time defined $\nabla$ naturally on all degrees. In the nondegenerate case we know from Lemma~\ref{symanti} that it is unique for $C=T=0$ and  in this case we also have a picture of $L_\delta$ as Lie bracket of vector fields, so in this case we have a formula for $\nabla$ that depends only on the metric,  akin to the familiar Koszul formula. One also has torsion freeness and metric-compatibility on all degrees and in a suitable sense.

\subsection{Divergence operator}

We conclude by supplying the partial inverse to Theorem~\ref{levi}. 

\begin{proposition}\label{div} Let $\Omega(A)$ be of classical type and $(\ ,\ )$ the inverse of a symmetric metric $g$ with metric compatible and torsion free covariant derivative $\nabla$ extending as a derivation to $\Omega$. Then the `divergence operator' 
\[ \delta_\nabla=\inter_{g\uo}\nabla_{g\ut}\]
is also of classical type and application of Theorem~\ref{levi} recovers $(\ ,\ ),\nabla$. 
\end{proposition}
\proof Here $\delta=\delta_\nabla$ as stated is well-defined as $\nabla_\omega$ is tensorial in $\omega$. Clearly 
\[ \delta(a\omega)=\inter_{g\uo}((g\ut,\extd a)\omega+a\nabla_{g\ut}\omega)=a\delta\omega+\inter_{g\uo}(\omega)(g\ut,\extd a)=a\delta\omega+\inter_{\extd a}\omega\] 
applies with the interior product provided by the given metric inner product, which is indeed a graded-derivation and symmetric. Hence Lemma~\ref{deriv} applies and $\delta$ is of classical type provided we can show that $\delta^2$ is an $A$-module map. As in the proof of Lemma~\ref{deriv}  (worked in reverse), this amounts to showing that $\inter_{\extd a}$ and $\delta$ anticommute for all $a\in A$. We do this in two steps. First, we observe that for any $\omega,\eta\in \Omega^1$,
\begin{equation}\label{comrel}[ \inter_{\eta},\nabla_\omega]=-\inter_{\nabla_\omega\eta}\end{equation}
holds as operations on $\Omega$. Indeed, using the Leibniz property of $\nabla_\omega$ and the graded-Leibniz property of interior one can check that the left hand side of (\ref{comrel}) is a degree -1 graded-derivation. The right hand side is also a graded-derivation and (\ref{comrel}) holds in degree 1 by metric compatibility (and is trivial in degree 0). Here for $\omega,\eta\in\Omega^1$,
\[ [\inter_{\eta},\nabla_\omega]\zeta=(\eta,\nabla_\omega\zeta)-(\omega,\extd(\eta,\zeta))=-(\zeta,\nabla_\omega\eta)\]
by metric compatibility (\ref{metriccomp}). Next, using (\ref{comrel}), we compute for all $a\in A, \omega\in \Omega$,
\[ (\delta\inter_{\extd a}+\inter_{\extd a}\delta)(\omega)=\inter_{g\uo}\nabla_{g\ut}\inter_{\extd a}\omega-\inter_{g\uo}\inter_{\extd a}\nabla_{g\ut}\omega =\inter_{g\uo}\inter_{\nabla_{g\ut}\extd a}(\omega)=0\]
since $g\uo\nabla_{g\ut}(\extd a)=0$ as an expression of zero torsion. Indeed, if $T=0$ then
\[\inter_\eta\inter_\omega( g\uo\nabla_{g\ut}\zeta)=\inter_\eta(\nabla_\omega\zeta-g\uo\inter_\omega\nabla_{g\ut}\zeta)=\inter_\eta\nabla_\omega\zeta-\inter_\omega\nabla_\eta\zeta=\inter_\eta\inter_\omega\extd\zeta.\]
 This concludes our proof that $(\Omega(A),\delta_\nabla)$ is of classical type. Now consider the covariant derivative defined by Theorem~\ref{levi}. It clearly coincides with the given $\nabla$ on degree 1 since both are metric compatible and torsion free and $(\ ,\ )$ is non-degenerate (see Lemma~\ref{symanti}).  Both covariant derivatives are derivations, in the case of the one in Theorem~\ref{levi} by Lemma~\ref{Schouten} and since  $\CL_\omega+\extd\omega\perp$ is a derivation, hence the two covariant derivatives coincide in all degrees.  \endproof

This implies in particular that every invertible metric (and associated covariant derivative) is in the image of the construction of Theorem~\ref{levi} for some choice of $\delta$. The same result applies more generally to $\nabla=\nabla\uo\tens\nabla\ut$ an algebraic connection and $(\ ,\ )$ nondegenerate, where we take $\delta_\nabla=\inter_{\nabla\uo\eta}(\nabla\ut\eta)$. We have focussed on the invertible case as the algebraic case involves more radically unfamiliar notation. 

Finally, for a fixed DGA $\Omega(A)$ of classical type let $\CD$ be the set of $\delta$ that are also of classical type (i.e., such that $(\Omega(A),\delta)$ is of classical type) and let   $\CD^\times$ be the subset where in addition the associated $(\ ,\ )$ is invertible. We consider Lemma~\ref{deriv} and Theorem~\ref{levi} more formally as providing a map $\Theta$ that associates to $\delta\in\CD$ a pair $\Theta(\delta)=((\ ,\ ),\nabla)$ consisting of a metric inner product $(\ ,\ )$ and a metric compatible torsion free covariant derivative $\nabla$ on $\Omega$, not necessarily a derivation. We denote by $\CI$ the set of vector fields, i.e. tensorial maps $v:\Omega^1\to A$ extended as degree -1 graded-derivations $\lfloor_v$. We call these `interior coderivations' and note that they form a vector space. 

\begin{proposition} \label{bijection} Let $\Omega(A)$ be of classical type and $\Theta$ the map in Theorem~\ref{levi}. Then $\CI$ acts on $\CD$ and $\CD^\times$ by addition and $\Theta$ becomes injective on the quotient $\CD/\CI$. Moreover, $\Theta$ restricts to a bijection from $\CD^\times/\CI$ to the set of invertible metrics and their associated covariant derivatives.  \end{proposition} 
\proof Let $\Omega(A)$ be of classical type and $\delta\in\CD$, and let $\lfloor_v$ an inner coderivation along $v:\Omega^1\to A$. Note that $\lfloor_v$ anticommutes with $\inter_\omega$ for all $\omega\in\Omega^1$. As $\lfloor_v$ is a degree -1 graded derivation, $\delta'=\delta+\lfloor_v$ has the same Leibnizator as $\delta$. Hence $\delta'$ is regular with the same interior products and such that $L_\delta(\omega,a)=L_{\delta'}(\omega,a)=\inter_{\extd a}$, so $\delta'\in\CD$ and it has the same metric by Lemma~\ref{deriv}. Here
\[ (\delta\lfloor_v+\lfloor_v\delta)(a\omega)=\delta(a\lfloor_v\omega)+\lfloor_v(a\delta\omega+\inter_{\extd a}\omega)=a(\delta\lfloor_v+\lfloor_v\delta)\omega+\{\inter_{\extd a},\lfloor_v\}\omega=a(\delta\lfloor_v+\lfloor_v\delta)\omega\]
so  $\delta'^2$ is a module map if $\delta^2$ is. $L_\delta=L_{\delta'}$ also means that the associated covariant derivatives in Theorem~\ref{levi} have the same first term and they have the same remaining terms as these depend only on $\extd,\inter_\omega$. Hence the covariant derivatives are the same in all degrees. We see that $\CI$ acts on $\CD$ and that $\Theta$ descends to the quotient. Conversely, suppose $\delta,\delta'\in \CD$ and $\Theta(\delta)=\Theta(\delta')$, i.e. they lead to the same metric inner product and covariant derivative. Then $\delta,\delta'$ have the same Leibnizator if one argument is in degree 0, as this is the interior product. More generally, as they result in the same covariant derivative on all degrees in Theorem~\ref{levi} we conclude that $L_\delta(\omega,\eta)=L_{\delta'}(\omega,\eta)$ for all $\omega\in \Omega^1$ and all $\eta\in \Omega$. 
Next we recall the tautological identity (\ref{L1}) for $L_\delta$, and the same for $L_{\delta'}$.  It follows by induction on the degree of the first argument that  $L_\delta=L_{\delta'}$ in all degrees. Hence $\delta'-\delta$ is a degree -1 graded-derivation. A  degree -1 graded-derivation is determined by its value on degree 1 as a map $\Omega^1\to A$ and hence takes the form of an interior coderivation (and necessarily anticommutes with any interior products). We see that  $\Theta$ is injective on the quotient space, i.e. the orbits of $\CI$ are precisely its fibres. As the action of $\CI$ does not change the associated metric, it acts on $\CD^\times$. By Proposition~\ref{div} every invertible metric (and associated metric compatible torsion free covariant derivative on $\Omega$ acting as derivations) is in the image of $\Theta$ and the latter now becomes a bijection with inverse provided by $\delta_\nabla$. \endproof

\begin{remark}\label{classical} Let $(M,g)$ be a (say) smooth orientable Riemannian manifold and $\star$ the Hodge duality defined by $(\omega,\eta){\rm Vol}=\omega\star\eta$ for the metric inner product $(\ ,\ )$ extended to forms of the same degree. Define $\delta$ by $\delta(\omega)=(-1)^{|\omega|+1}\star^{-1}\extd \star(\omega)$. This implies
\[ \delta(a\omega)=(-1)^{|\omega|+1}\star^{-1}\extd\star(a\omega)=a\delta\omega+(-1)^{|\omega|+1}\star^{-1}((\extd a)\star\omega)\]
so that
\[ \inter_{\omega}\eta=(-1)^{|\eta|-1}\star^{-1}(\omega (\eta^\star)),\quad\forall\omega\in \Omega^1,\ \eta\in\Omega,\]
which is a known formula for the interior product and known to provide a left-derivation, so Lemma~\ref{deriv} applies. These conventions mean that $\delta$ is adjoint to $-\extd$ in the sense of Hodge theory, which is a known convention though not necessarily the most popular one. In this convention the Hodge Laplacian and the Laplace-Beltrami operators coincide in degree 0 rather than with a minus sign. Note also that when $\inter_\omega$ is extended to $\omega\in \Omega$ as discussed above, it again has the form of $\star^{-1}(\omega(\ )^\star)$ but with a more complicated sign factor and that in particular $\inter_\omega\eta=(-1)^{m(m-1)\over 2}(\omega,\eta)$ for $\omega,\eta\in \Omega^m$. Clearly Theorem~\ref{levi} again recovers the Levi-Civita covariant derivative and provides a  formula for it in terms of the Lie derivative, interior product and the failure of $\delta$ to be a graded-derivation. Because the same applies to $\delta_\nabla$ in Proposition~\ref{div}, we conclude that the two coincide possibly up to  interior product along a vector field. In fact there is no such vector field as $\delta=\delta_\nabla$  is equivalent, given the above, to $\extd\omega^\star=g\uo(\nabla_{g\ut}\omega)^\star$ for all $\omega\in\Omega$. From the formula in Theorem~\ref{levi} it is easy to see that $\star$ commutes with $\nabla$ (there are also other easy ways to see this) so we require $\extd\omega=g\uo\nabla_{g\ut}\omega$ for all $\omega\in \Omega$. But on degree 1 this is just the content of zero torsion (see the proof of Propsition~\ref{div}) and hence holds in all degrees by derivation properties of both sides. Another modest application is to the Leibniz property of the Hodge-Laplacian $\Delta=\extd\delta+\delta\extd$. We note the tautological identity 
\begin{equation}\label{L2} L_\Delta(\omega,\eta)=\extd L_\delta(\omega,\eta)+L_\delta(\extd\omega,\eta)+(-1)^{|\omega|}L_\delta(\omega,\extd\eta),\quad\forall\omega,\eta\in\Omega\end{equation}
valid for any degree -1 linear map $\delta$ on any DGA on writing out the definitions of all the terms, and as observed in \cite{Coll} in the present graded-commutative context. A special case is $\Delta(a\omega)=(\Delta a)\omega+a\Delta(\omega)+L_\delta(\extd a,\omega)+\CL_{\extd a}\omega$ for all $a\in A$ and $\omega\in \Omega$. By Theorem~\ref{levi} the last two terms are $2\nabla_{\extd a}\omega$, giving a 2nd order Leibniz rule normally proven by other means.  \end{remark}

\subsection{Ricci tensor}\label{secricci}

Here we give a more substantial application of  our new formula for  the Levi-Civita covariant derivative in Theorem~\ref{levi}. We suppose that $(\Omega(A),\delta)$ is of classical type with $(\ ,\ )$ invertible, with inverse $g$, and we let $\Delta=\extd\delta+\delta\extd$. We assume that $\delta=\delta_\nabla$ as we know can be arranged in the classical setting, see Remark~\ref{classical}. 

\begin{lemma}\label{second} Let $B$ be a degree 0 linear map on $\Omega$ such that $L_B(a,\omega)=2\nabla_{\extd a}\omega$ for all $a\in A,\ \omega\in\Omega$. Then $B$ extends canonically to tensor products as 
\[ B(\omega\tens_A\eta)=B\omega\tens_A\eta+\omega\tens_AB\eta+2\nabla_{g\uo}\omega\tens_A \nabla_{g\ut}\eta\]
for all $\omega,\eta\in \Omega$.
\end{lemma}
\proof  First note that the construction is depends tensorially on $g$, i.e. only on $g\in \Omega^1\tens_A\Omega^1$ so it is well-defined.  With $\tens=\tens_A$,
\begin{eqnarray*} 
 B(\omega a \tens \eta)&=&B\omega\tens a\eta+2\nabla_{\extd a}\omega\tens\eta+\omega\tens (Ba)\eta+\omega\tens a B \eta+2\nabla_{g\uo}a\omega\tens \nabla_{g\ut}\eta\\
 &=&B\omega\tens a\eta+2\nabla_{\extd a}\omega\tens\eta+\omega\tens B(a\eta)-\omega\tens 2\nabla_{\extd a}\eta\\
 &&+\nabla_{g\uo}\omega\tens \nabla_{g\ut}\eta+2\omega\tens \nabla_{\extd a}\eta\\
 &=&B\omega\tens a\eta+2\nabla_{\extd a}\omega\tens \eta+\omega\tens  B(a\eta)+2\nabla_{g\uo}\omega\tens \nabla_{g\ut}a\eta-2 \nabla_{\extd a}\omega\tens\eta\\
 &=& B(\omega  \tens a \eta)
 \end{eqnarray*}
 so that the construction stated descends to a map on $\Omega\tens_A\Omega$. \endproof
 
The Leibnizator here is characteristic of a covariant second-order operator and holds for the Hodge-Laplacian $\Delta$  by Remark~\ref{classical} (the explanation is valid for any
$(\Omega(A),\delta)$ of classical type). 

\begin{proposition}  The `Laplace-Beltrami' operator (\ref{LB}) obeys
\[ L_{\Delta_{LB}}(\omega,\eta)=2(\nabla_{g\uo}\omega)\nabla_{g\ut}\eta\]
for all $\omega,\eta\in \Omega$. In particular, it is of the type in Lemma~\ref{second}. Moreover,

(1) $W:=\Delta_{LB}-\Delta$ is tensorial (an $A$-module map) and zero in degree 0.

(2) $\Delta_{LB}(g)=0$.
\end{proposition}
\proof We have already explained the definition of $\Delta_{LB}$ in the Preliminaries in full generality. In our case we compute
\begin{eqnarray*}
\Delta_{LB}(\omega\eta)&=&\nabla_{g\uo}\left((\nabla_{g\ut}\omega)\eta+\omega\nabla_{g\ut}\eta\right)-(\nabla_{\nabla_{g\uo}g\ut}\omega)\eta-\omega\nabla_{\nabla_{g\uo}g\ut}\eta\\
&=&(\Delta_{LB}\omega)\eta+\omega\Delta_{LB}\eta+(\nabla_{g\uo}\omega)\nabla_{g\ut}\eta+ (\nabla_{g\ut}\omega)\nabla_{g\uo}\eta.
\end{eqnarray*}
The last two terms are the same since they depend tensorially on $g$ and hence we can use its symmetry. As a special case, we see that  $L_{\Delta_{LB}}(a,\omega)=2\nabla_{\extd a}\omega$ for all  $a\in A,\ \omega\in \Omega$, so Lemma~\ref{second} applies. 

Next, (1) In degree 0, $\Delta a=\delta\extd a=\inter_{g\uo}\nabla_{g\ut}\extd a=\nabla_{g\ut}\inter_{g\uo}\extd a-\inter_{\nabla_{g\ut}g\uo}\extd a=\inter_{g\ut}\extd\inter_{g\uo}\extd-\inter_{\nabla_{g\ut}g\uo}\extd a=\Delta_{LB}a$ for all $a\in A$, using metric compatibility and symmetry of the metric. Hence $W$ is zero on degree 0. Since both $\Delta,\Delta_{LB}$ have the same Leibnizator when one argument is in $A$, we then have $W(a\omega)=W(a)\omega+a W(\omega)=aW(\omega)$. 

(2) We evaluate  half of the desired expression against $\omega\tens\eta$,
\begin{eqnarray*} &&\kern-20pt(\omega,\nabla_{g\uo}\nabla_{g\ut}g'\uo)(\eta,g'\ut)+(\omega,\nabla_{g\uo}g'\uo)(\eta,\nabla_{g\ut}g'\ut)\\
&=&(\omega,\nabla_{g\uo}((\eta,g'\ut)\nabla_{g\ut}g'\uo))-(g\uo,\extd(\eta,g'\ut))(\omega,\nabla_{g\ut}g'\uo)\\
&&+(\omega,\nabla_{g\uo}g'\uo)(-(\nabla_{g\ut}\eta,g'\ut)+(g\ut,\extd(\eta,g'\ut))\\
&=&(\omega,(\nabla_{g\uo}((\eta,g'\ut)\nabla_{g\ut}g'\uo))-(\omega,\nabla_{g\uo}g'\uo)(\nabla_{g\ut}\eta,g'\ut)\\
&=&-(\omega,\nabla_{g\uo}((g\ut,\extd(\eta,g'\ut))g'\uo))+(\omega,g'\uo)(g\uo,\extd (\nabla_{g\ut}\eta,g'\ut))\\
&=&-(\omega,\nabla_{\extd(\eta,g'\ut)} g'\uo)+(\omega,g'\uo)(g\uo,\extd (\nabla_{g\ut}\eta,g'\ut))-(\omega,g'\uo)(g\uo,\extd(g\ut,\extd(\eta,g'\ut)))\\
&=&-(\omega,\nabla_{\extd(\eta,g'\ut)} g'\uo)+(g\uo,\extd(\nabla_{g\ut}\eta,\omega))-(\nabla_{\extd(\omega,g'\uo)}\eta,g'\ut)\\
&&-(g\uo,\extd (g\ut,(\omega,g'\uo)\extd(\eta,g'\ut)))+(\extd(\omega,g'\uo),\extd(\eta,g'\ut)))
\end{eqnarray*}
where for the first equality we moved a scalar factor inside a covariant derivative and compensated and we used metric compatibility
to move over to an action on $\eta$. We repeat the first process so as to be able to cancel a metric with its inverse, and repeat this principle.  A similar calculation for the other half gives
\[ -(\eta,\nabla_{\extd(\omega,g'\uo)} g'\ut)+(g\uo,\extd(\nabla_{g\ut}\omega,\eta))-(\nabla_{\extd(\eta,g'\ut)}\omega,g'\uo)\]
\[ -(g\uo,\extd (g\ut,(\eta,g'\ut)\extd(\omega,g'\uo)))+(\extd(\eta,g'\ut),\extd(\omega,g'\uo)))\]
Adding these together using metric compatibility and the Leibniz rule for $\extd$ gives zero.
\endproof

Finally, we prove the relationship with the Ricci tensor. We define the Ricci map by $\widetilde{\rm Ricci}(\omega)=R(\omega,g\uo)g\ut$ for all $\omega\in \Omega^1$, where $R$ is the Riemann curvature. This is well-defined and tensorial by the tensoriality of Riemann. The Ricci tensor itself is then defined by ${\rm Ricci}=g\uo\tens_A \widetilde{\rm Ricci}(g\ut)\in \Omega^1\tens_A\Omega^1$. 

\begin{corollary} At least in the case of a classical Riemannian manifold $(M,g)$,  ${\rm Ricci}=-{1\over 2}\Delta (g)$.
\end{corollary}
\proof At least in the classical case one knows that $W=\widetilde{{\rm Ricci}}$ on $\Omega^1$ and also that ${\rm Ricci}$ is symmetric. Then 
\[ \Delta(g)=\Delta(g\uo)\tens g\ut+g\uo\tens\Delta(g\ut)+2\nabla_{g\uo}g'\uo\tens \nabla_{g\ut}g'\ut\]
\[ =-W(g\uo)\tens g\ut-g\uo\tens W(g\ut)+\Delta_{LB}(g)=-2{\rm Ricci}\]
using the symmetry proven.  This then provides the stated formula for the Ricci tensor. \endproof

We expect the same result for all $(\Omega(A),\delta)$ of classical type, using the methods as above. The required symmetry of Ricci is straightforward to prove but the calculation of $W$ using our particular methods appears to be more tedious.

\section{Riemannian structures induced by central extensions}

In this section we see how a metric and covariant derivative arise naturally from an extension problem in noncommutative geometry, including how the datum $\delta$ in Section~2 arises naturally. 

\subsection{Central extensions of DGAs} 

We first formulate the required notion of a `central extension' of a general DGA $\Omega(A)$ in degree 1 by the algebra $\Omega_{\theta'}=k[\theta']/\<\theta'^2\>$  viewed as a trivial DGA with $\theta'$ of degree 1 and $\extd\theta'=0$.

\begin{definition}\label{centralext} By central extension of a DGA $\Omega(A)$ we mean a DGA $\tilde{\Omega}(A)$ such that 
\[  \tilde{\Omega}(A)= \Omega_{\theta'}\tens\Omega(A)\] as a vector space and 
\[ 0\to \Omega_{\theta'}\to \tilde{\Omega}(A)\to \Omega(A)\to 0\]
as maps of DGA's, where the maps come from the canonical inclusion  in the tensor product and by setting $\theta'=0$. We also require that $\Omega_{\theta'}$ here is graded-central,
\[ \theta'\omega=(-1)^{|\omega|}\omega\theta'\]
in $\tilde{\Omega}(A)$. A morphism of extensions $\Phi:\tilde{\Omega}(A)\to \tilde{\Omega'}(A)$ means a map of DGA's such that
\[ \Omega_{\theta'} \begin{array}{ccc}  & \tilde{\Omega}(A) & \\ {\buildrel \nearrow\over\searrow} & \quad\downarrow\Phi & {\buildrel \searrow\over\nearrow} \\ & \tilde{\Omega'}(A) &\end{array}\Omega(A),\quad \Phi(\theta')=\theta',\quad\Phi(\omega)=\omega - {\lambda\over 2}\theta'\delta(\omega)\] 
By a (left) {\em cleft} central extension we mean a central extension where the canonical linear inclusion of  $\Omega(A)$ coming from the tensor product form is a left $A$-module map. 
\end{definition}
Clearly the exterior derivative and product of $\tilde{\Omega}(A)$ must necessarily have the form
\[ \omega\cdot\eta=\omega\eta-{\lambda\over 2}\theta'\lo\omega,\eta\lc,\quad \extd_\cdot\omega=\extd \omega - {\lambda\over 2}\theta' \Delta\omega,\quad\omega,\eta\in\Omega(A)\]
for a bilinear map $\lo\ ,\ \lc$ of degree -1 and a linear map $\Delta$ of degree 0. This form is necessary since $\theta'$ has degree 1. The $\lambda/2$ is a parameter which we insert here in the normalisations of the maps as it may be relevant to a future deformation analysis,  but for our purposes we think of it as a non-zero element of the ground field and can set it to 1. The extension is cleft precisely when $\lo a,\ \lc=0$ for all $a\in A$. 

\begin{proposition}\label{Omegaext} Let $\Omega(A)$ be a DGA on an algebra $A$. Degree 0,-1 maps $\Delta:\Omega(A)\to \Omega(A)$ and $\lo\ ,\ \lc:\Omega(A)\tens\Omega(A)\to \Omega(A)$ respectively define a central extension $\tilde{\Omega}(A;\Delta,\lo\ ,\ \lc)$ iff 
\begin{equation}\label{c1} \lo\omega\eta,\zeta\lc+\lo\omega,\eta\lc\zeta=\lo\omega,\eta\zeta\lc+(-1)^{|\omega|}\omega\lo\eta,\zeta\lc.\end{equation}
\begin{equation}\label{c2} L_\Delta(\omega,\eta)=\extd\lo\omega,\eta\lc+\lo\extd\omega,\eta\lc+(-1)^{|\omega|}\lo\omega,\extd\eta\lc
\end{equation}
for all $\omega,\eta,\zeta\in \Omega(A)$, and $[\Delta,\extd]=0$. 
\end{proposition}
\proof For associativity we compute 
\begin{eqnarray*} (\omega\cdot\eta)\cdot\zeta&=&\left(\omega\eta+{\lambda\over 2}(-1)^{|\omega|+|\eta|}\lo\omega,\eta\lc\theta'\right)\cdot\zeta\\
&=&\omega\eta\zeta+{\lambda\over 2}(-1)^{|\omega|+|\eta|+|\zeta|}\lo\omega,\eta\lc\zeta\theta'+{\lambda\over 2}(-1)^{|\omega|+|\eta|+|\zeta|}\lo\omega\eta,\zeta\lc\theta'\\
\omega\cdot(\eta\cdot\zeta)&=&\omega\cdot\left(\eta\zeta+{\lambda\over 2}(-1)^{|\eta|+|\zeta|}\lo\eta,\zeta\lc\theta'\right)\\
&=&\omega\eta\zeta+{\lambda\over 2}(-1)^{|\eta|+|\zeta|}\omega\lo\eta,\zeta\lc\theta'+ {\lambda\over 2}(-1)^{|\omega|+|\eta|+|\zeta|}\lo\omega,\eta\zeta\lc\theta'
\end{eqnarray*}
Comparing, we see that we need (\ref{c1}).  Next, for the Leibniz rule we compute
\begin{eqnarray*}\extd_\cdot(\omega\cdot\eta)&=&\extd_\cdot\left( \omega\eta+{\lambda\over 2}(-1)^{|\omega|+|\eta|}\lo\omega,\eta\lc\theta' \right) \\
&=&\extd(\omega\eta)-{\lambda\over 2}(-1)^{|\omega|+|\eta|}(\Delta(\omega\eta))\theta'+{\lambda\over 2}(-1)^{|\omega|+|\eta|}(\extd \lo\omega,\eta\lc)\theta'\\
(\extd_\cdot\omega)\cdot\eta+(-1)^{|\omega|}\omega\cdot\extd_\cdot\eta&=&\left(\extd\omega-{\lambda\over 2}(-1)^{|\omega|}(\Delta\omega)\theta'\right)\cdot\eta+(-1)^{|\omega|}\omega\cdot\left(\extd \eta-{\lambda\over 2}(-1)^{|\eta|}(\Delta\eta)\theta'\right)\\
&=&(\extd\omega)\eta+(-1)^{|\omega|}\omega\extd\eta+{\lambda\over 2}(-1)^{|\omega|+1+|\eta|}\lo\extd\omega,\eta\lc\theta'+{\lambda\over 2}(-1)^{|\eta|+1}\lo\omega,\extd\eta\lc\theta'\\
&&-{\lambda\over 2}(-1)^{|\omega|+|\eta|}(\Delta\omega)\eta\theta'- {\lambda\over 2}(-1)^{|\eta|+|\omega|}\omega(\Delta\eta)\theta'
\end{eqnarray*}
where we used $\theta'^2=0$ and $\extd_\cdot\theta'=0$. Comparing, we see that we need (\ref{c2}). We also need \[\extd_\cdot\extd_\cdot\omega=\extd_\cdot(\extd \omega-{\lambda\over 2}(-1)^{|\omega|}(\Delta\omega)\theta')= \extd^2\omega-{\lambda\over 2}(-1)^{|\omega|+1}(\Delta\extd\omega)\theta'-{\lambda\over 2}(-1)^{|\omega|}(\extd\Delta\omega)\theta'=0\]
which requires $[\Delta,\extd]=0$. \endproof

We refer to the  pair $(\Delta,\lo\ ,\ \lc)$ obeying (\ref{c1})-(\ref{c2}) and $[\Delta,\extd]=0$ as a  {\em 2-cocycle} on the DGA  in analogy with the way that central extensions of groups are defined by 2-cocycles. To complete this picture:

\begin{lemma}\label{equivsols} Two cocycles $(\Delta, \lo\ ,\ \lc)$ and $(\Delta', \lo\ ,\ \lc')$  give isomorphic central extensions iff
\[ \Delta'= \Delta + \extd\delta+\delta\extd,\quad  \lo\omega,\eta\lc'=\lo\omega,\eta\lc+L_\delta(\omega,\eta)\]
for some degree -1 linear map $\delta: \Omega(A)\to \Omega(A)$. 
\end{lemma} 
\proof Any map  $\Phi:\tilde{\Omega}(A;\Delta,\lo\ ,\ \lc)\to \tilde{\Omega}(A;\Delta',\lo\ ,\ \lc')$ that commutes with the canonical inclusions and projections (a morphism of extensions) must have the form:
\[ \omega\mapsto \omega-{\lambda\over 2}\theta'\, \delta\omega,\quad \theta'\mapsto \theta',\quad\forall \omega\in\Omega(A),\]
for some degree -1 map $\delta$. One may easily verify from the definitions  that this is an isomorphism iff the difference between the $\cdot$ products and $\extd_\cdot$ in the two cases have the form stated. \endproof

We refer a cocycle of the form
\begin{equation}\label{itocalc} \Delta=\extd\delta+\delta \extd,\quad \lo\omega,\eta\lc=L_\delta(\omega,\eta),\quad\forall \omega,\eta\in \Omega.\end{equation}
associated to any degree -1 linear map $\delta$ as its {\em coboundary} and in this case  Lemma~\ref{equivsols} says that central extensions up to isomorphism are classified by cocycles $(\Delta,\lo\ ,\ \lc)$ up to such coboundaries, i.e. by a form of 2-cohomology. We have not placed the cohomology here into a general context but this is parallel to  central extensions of a group being classified by its 2-cohomology. 

We have already observed in Section~2 that (\ref{itocalc}) tautologically solves (\ref{c1})-(\ref{c2}) for any degree -1 linear map $\delta$, cf \cite{Coll} in the graded-commutative case.  This `homologically trivial' case is not interesting from our point of view of `quantisation' as it does not change the DGA but it can still be of interest, see \cite{Beggs} where a similar DGA in this case has been noted and applied to the understanding of the It\^o stochastic differential cf. \cite{MH}.  

We now restrict this degree of equivalence by focussing on cleft extensions, and will then see how Riemannian geometry `emerges' from this restricted extension problem.

\begin{lemma}\label{uniquesol} Let $\Omega(A)$ be a standard DGA. In a cleft central extension the $\lo\ ,\ \lc$ part of the cocycle is uniquely determined by the $\Delta$ part.
\end{lemma}
\proof Here we are supposing that $\lo a,\  \lc=0$ for all $a\in A$. In this case we have
\begin{equation}\label{l0} L_\Delta(a,\eta)=\lo \extd a,\eta\lc,\quad\forall a\in A,\quad\eta\in \Omega(A)\end{equation}
is a special case of (\ref{c2}). Next, we specialise (\ref{c1}) to $a\in A$, $\omega,\eta,\zeta\in\Omega(A)$ as
\begin{eqnarray}\label{l1} \lo a\eta,\zeta\lc&=&a\lo \eta,\zeta\lc\\\
\label{l2} \lo \omega a,\zeta\lc +\lo \omega,a\lc \zeta&=&\lo \omega, a\zeta\lc \\
\label{l3} \lo \omega\eta,a\lc +\lo \omega,\eta\lc a&=&\lo\omega,\eta a\lc + (-1)^{|\omega|}\omega\lo \eta,a\lc
\end{eqnarray}
We now proceed as follows. From $L_\Delta$ we define $\lo\extd a,\eta\lc$ for all $\eta$. By (\ref{l1}) and the assumption that the DGA is surjective (the standard case) we see that we have defined $\lo\omega,\eta\lc$ for all $\omega\in \Omega^1$ and all $\eta$. Now suppose for a fixed $\zeta\in \Omega(A)$ that $\lo \eta ,\zeta\lc $ has  been defined up to $\eta$ of some degree. Then  (\ref{c1}) defines $\lo \omega\eta,\zeta\lc$ for any $\omega$. In this way, assuming $\Omega(A)$ is generated by degree 0 and 1, we have defined $\lo\eta,\zeta\lc$ for all $\eta$. The initial case for the induction where $\eta$ has degree 1 was already specified for any $\zeta$ earlier in the construction.  \endproof

For completeness we also list a remaining special case of (\ref{c2}),
\begin{equation}\label{l00} L_\Delta(\omega,a)=\extd\lo \omega,a\lc + \lo\extd\omega,a\lc+ (-1)^{|\omega|}\lo\omega,\extd a\lc,\quad\forall\omega\in\Omega,\ a\in A\end{equation}
which will need  later.

\subsection{Bimodule covariant derivatives associated to cleft extensions}

\begin{definition} \label{cleftreg} We say that a cleft extension $(\Delta,\lo\ ,\ \lc)$ on a standard DGA $\Omega(A)$ is $n$-{\em regular} if
\[ \cj_\omega(a\extd b)={1\over 2}\lo \omega a,b\lc,\quad \forall \omega\in \Omega,\ a,b\in A\]
is  a well-defined degree -1 map $\cj:\Omega^i \tens_A\Omega^1\to \Omega^{i-1}$ for $i\le n$. We say that the cleft extension is regular is it is regular for all degrees. We refer to $\cj$ as `interior product' and its restriction $(\ ,\ ):\Omega^1\tens_A\Omega^1\to A$  to degree 1 as `metric'. 
\end{definition}

The following is stated for regular extensions but  we need only 1-regularity for the metric and a covariant derivative to be defined and 2-regularity for this to be a bimodule covariant derivative acting on degree 1 (which is the main case of interest for Riemannian geometry). 

\begin{proposition}\label{cleftcon1} If $(\Delta,\lo\ ,\ \lc)$ is a regular cleft extension on a standard DGA $\Omega(A)$ then $\cj$ is a bimodule map and
\[ \nabla_\omega\eta={1\over 2} \lo \omega,\eta\lc,\quad\forall \omega\in \Omega^1,\quad \eta\in\Omega\]
is a bimodule covariant derivative on $\Omega$ with respect to $(\ ,\ )$.  Here
\[ \sigma:\Omega^1\tens_A\Omega\tens_A\Omega^1\to \Omega,\quad \sigma_\omega(\eta\tens_A\zeta)=\cj_{\omega\eta}(\zeta)+\omega\cj_\eta(\zeta),\quad\forall \omega,\zeta\in \Omega^1,\ \eta\in \Omega.\]
Moreover, if $(\ ,\ )$ is invertible with metric $g$ then $\nabla_\omega$  acts on tensor products,
\begin{equation}\label{covtens} \nabla_\omega(\eta\tens_A\zeta)=\nabla_\omega\eta\tens_A\zeta+ \sigma_\omega(\eta\tens_A g\uo )\tens_A\nabla_{g\ut  }\zeta,\quad\forall \eta,\zeta\in \Omega.\end{equation}
\end{proposition}
\proof On 0-forms (\ref{l1}) tells us that $\cj_{a\omega}=a\cj_\omega$ and (\ref{l2}) combined with the Leibniz rule tells us that $\cj_\omega((\extd a)b)=\cj_\omega(\extd a)b$ so that $\cj$ becomes a bimodule map. On 1-forms (\ref{l1}) tells us that $\nabla_{a\omega}=a\nabla_\omega$. Meanwhile (\ref{l2}) tells us that $\nabla_\omega$ 
is a covariant derivative in the sense (\ref{leftcov}) given the definition of $(\ ,\ )$. Given this, (\ref{l3})  tells us that we have a bimodule covariant derivative in the sense (\ref{leftbimod}) provided
we define $\sigma$ as stated. That $\sigma$ is a bimodule map is immediate from the properties of $\cj$. If $(\ ,\ )$ is invertible then a bimodule covariant derivative 
induces an algebraic bimodule connection and these act on tensor products as in (\ref{covtens}).  \endproof

One might expect in the invertible case that the above action on tensor products is compatible with $\wedge$, which is to say:
\begin{equation}\label{sigmaderiv} \nabla_\omega(\eta\zeta)=(\nabla_\omega\eta)\zeta+\sigma_\omega(\eta\tens_A g\uo )\nabla_{g\ut  }\zeta,\quad\forall \eta,\zeta\in \Omega\end{equation}
and this can be the case but does not appear always to be true. However, there is always a different kind of generalised Leibniz rule:

\begin{proposition}\label{cleftcon2} Given a regular cleft extension, the general $\nabla_\omega={1\over 2}\lo\omega,\ \lc$ for all $\omega\in \Omega$ is a left-covariant derivative in the sense
\[ \nabla_{a\omega}=a\nabla_\omega,\quad \nabla_\omega(a\eta)=\nabla_{\omega a}\eta+ \cj_\omega(\extd a)\eta\]
and a bimodule covariant derivative in the sense (\ref{leftbimod}) but now with
\[ \sigma:\Omega\tens_A\Omega\tens_A\Omega^1\to \Omega,\quad \sigma_\omega(\eta\tens_A\zeta)= \cj_{\omega\eta}(\zeta)- (-1)^{|\omega|}\omega\cj_\eta(\zeta).\]
Moreover,
\begin{equation}\label{c1con} \nabla_\omega(\eta\zeta)=(\nabla_\omega\eta)\zeta- (-1)^{|\omega|}\omega\nabla_\eta\zeta+\nabla_{\omega\eta}\zeta,\quad\forall \omega,\eta,\zeta\in \Omega. \end{equation}
\begin{equation}\label{c2con} {1\over 2}L_\Delta(\omega,\eta)=(\extd\nabla_\omega+(-1)^{|\omega|}\nabla_\omega\extd)\eta+\nabla_{\extd\omega}\eta,\quad\forall \omega,\eta,\in \Omega. \end{equation}
\end{proposition}
\proof That we have a generalised bimodule covariant derivative follows exactly the same argument as the proof of Proposition~\ref{cleftcon1}, just now using more general forms in (\ref{l1})-(\ref{l3}). Moreover, the general (\ref{c1})  and (\ref{c2}) now become the two displayed equations (\ref{c1con}) and (\ref{c2con}) respectively. \endproof

Also observe that when $(\ ,\ )$ is invertible and noting that in this case the metric $g$ is necessarily central, it is easy to see that there is a potentially different higher-form bimodule covariant derivative
\[ \nabla'_\omega=\cj_\omega(g\uo )\nabla_{g\ut  },\quad \sigma'_\omega(\eta\tens_A\zeta)=\cj_\omega(g\uo )\sigma_{g\ut  }(\eta\tens_A\zeta)\]
which coincides with $\nabla_\omega$ for $\omega\in \Omega^1$.

\begin{proposition}\label{braleib} When $(\ ,\ )$ is invertible, the following are equivalent
\begin{enumerate}\item $(\nabla'_\omega,\sigma'_\omega)$ obey the braided-Leibniz rule (\ref{sigmaderiv}) for all $\omega\in\Omega^1$.
 \item $(\nabla_\omega,\sigma_\omega)$ obey the braided-Leibniz rule (\ref{sigmaderiv}) for all $\omega\in \Omega^1$.
 \item $\nabla_\omega=\nabla'_\omega$ for all $\omega\in \Omega$.
\end{enumerate}
In this case the braided-Leibniz rules also hold for all $\omega\in \Omega$ and 
\[\cj_{\omega\eta}=\cj_\omega(g\uo)(\cj_{g\ut\eta}+g\ut\cj_\eta)+(-1)^{|\omega|}\omega\cj_\eta,\qquad\forall \omega,\eta\in \Omega.\]
\end{proposition}
\proof That $(\nabla',\sigma')$ obeys (\ref{sigmaderiv}) is 
\[  \cj_\omega(g\uo)\nabla_{g\ut}(\eta\zeta)= \cj_\omega(g\uo)((\nabla_{g^2}\eta)\zeta+\sigma_{g\ut}(\eta\tens \bar g\uo)\nabla_{\bar g\ut}\zeta)\]
where $\bar g$ is another copy of $g$. Putting in the properties of $\nabla$ along 1-forms and cancelling, our condition is 
\[ \cj_\omega(g\uo)(g\ut\nabla_\eta\zeta+\nabla_{g\ut\eta}\zeta)= \cj_\omega(g\uo)(\cj_{g\ut\eta}(\bar g\uo)+g\ut\cj_\eta(\bar g\uo))\nabla_{\bar g\ut}\zeta\]
which holds for all $\omega\in\Omega$ (and all $\eta,\zeta$ understood) iff it holds for all $\omega\in \Omega^1$, where it reduces to the condition 
 \[ \nabla_{\omega\eta}-\cj_{\omega\eta}(g\uo)\nabla_{g\ut}=(-1)^{|\omega|}\omega(\nabla_\eta-\cj_\eta(g\uo)\nabla_{g\ut}),\quad\forall \omega\in\Omega^1,\quad \eta\in \Omega\]
 On the other hand, this is the condition that $(\nabla,\sigma)$ obeys (\ref{sigmaderiv}) given the form of $\sigma$ and moreover, if it holds for all $\omega\in\Omega^1$ then by induction on the degree of $\omega$, we conclude that $\nabla_\omega-\nabla'_\omega=0$ for all degrees of $\omega$ (given that this vanishes on $\omega$ of degree 1). In this case the last displayed condition vanishes identically for all degrees of $\omega$, so $(\nabla,\sigma)$ obeys $(\ref{sigmaderiv})$ for all degrees of $\omega$.  Finally, if $\nabla_\omega=\nabla'_\omega$ then $\sigma_\omega=\sigma'_\omega$  since these are uniquely determined, which is the stated condition on $\cj$. \endproof

 We also see that when the braided-Leibniz rule does hold, $\nabla_\omega$ along higher forms is given in terms of $\nabla$ along forms of lower degree and hence inductively in terms of the 1-form covariant derivative. Example~\ref{Z2} below is an instance where the braided-Leibniz rule holds and the above applies.

We conclude by studying some basic elements of the noncommutative geometry for this class of bimodule covariant derivatives.

\begin{proposition}\label{clefttorsion} Let $(\nabla_\omega,\sigma_\omega)$ be a bimodule 1-form covariant derivative as in Proposition~\ref{cleftcon1} and $(\ ,\ )$ be invertible. Then the algebraic torsion $T=g\uo \nabla_{g\ut  }-\extd:\Omega\to \Omega$ is a bimodule map (one says\cite{BegMa2} `torsion compatible') iff 
\[ g\uo g\ut   \cj_\omega(\zeta)+g\uo \cj_{g\ut  \omega}(\zeta)=(-1)^{|\omega|}\omega\zeta,\quad\forall\omega\in \Omega,\ \zeta\in \Omega^1.\]
If this holds and if the braided-Leibniz rule (\ref{sigmaderiv}) holds then $T$ is a derivation. 
\end{proposition}
\proof The torsion is already a left module map by the connection property and centrality of the metric. For the right module property we use the form of $\sigma$ to compute
\[ T(\omega a)=g\uo((\nabla_{g\ut}\omega)a+\sigma_{g\ut}(\omega\tens\extd a))-(\extd\omega)a-(-1)^{|\omega|}\omega\extd a\]
\[=(T\omega)a+g\uo\cj_{g\ut\omega}(\extd a)+g\uo g\ut\cj_\omega(\extd a)-(-1)^{|\omega|}\omega\extd a\]
so we require the condition stated. If this holds and (\ref{sigmaderiv}) holds then
\begin{eqnarray*} T(\omega\eta)&=&g\uo((\nabla_{g\ut}\omega)\eta+\sigma_{g\ut}(\omega\tens_A\bar g\uo)\nabla_{\bar g\ut}\eta)-\extd(\omega\eta)\\
&=&(T\omega)\eta-(-1)^{|\omega|}\omega\extd\eta+g\uo(\cj_{g\ut\omega}(\bar g\uo)+g\ut\cj_\omega(\bar g\uo)\nabla_{\bar g\ut}\\
&=& (T\omega)\eta-(-1)^{|\omega|}\omega\extd\eta+(-1)^{|\omega|}\omega g\uo\nabla_{g\ut}\eta=(T\omega)\eta+(-1)^{|\omega|}\omega T\eta\end{eqnarray*}
by the right-module condition. 
\endproof

\begin{proposition}\label{cleftricci}Let $(\nabla_\omega,\sigma_\omega)$ be a bimodule 1-form covariant derivative as in Proposition~\ref{cleftcon1} and $(\ ,\ )$ be invertible. The Laplace-Beltrami operator (\ref{LB}) has Leibnizator
\[ L_{\Delta_{LB}}(a,\omega)=\nabla_{2\extd a+\cj_{g\uo g\ut  }(\extd a)}\omega.\]
If $\wedge(g)=0$ (one says that $g$ is `quantum symmetric') then $\Delta_{LB}-\Delta$ is a left-module map and 
\[ {\rm Ricci}_\Delta:=g\uo \tens_A(\Delta_{LB}-\Delta)(g\ut  )\]
is well-defined. 
\end{proposition}
\proof For 1-form covariant derivative on a DGA $\Omega$ and $(\ ,\ )$ invertible, one has $\Delta_{LB}$ defined as explained in the Preliminaries. When we have a bimodule covariant derivative then we also have
\begin{eqnarray*}\Delta_{LB}(a\omega)&=&\nabla_{g\uo}\nabla_{g\ut}(a\omega)-\nabla_{\nabla_{g\uo}g\ut}(a\omega)\\
&=&\nabla_{g\uo}\nabla_{g\ut a}\omega-\nabla_{(\nabla_{g\uo}g\ut)a}\omega+\nabla_{g\uo}((\cj_{g\ut}(\extd a)\omega)-\cj_{\nabla_{g\uo}g\ut}(\extd a)\omega\\
&=&\nabla_{g\uo}\nabla_{g\ut a}\omega-\nabla_{\nabla_{g\uo}(g\ut a)}\omega+\nabla_{\sigma_{g\uo}(g\ut\tens_A\extd a)}\\
&&+\nabla_{g\uo \cj_{g\ut}(\extd a)}\omega+\cj_{g\uo}(\extd j_{g\ut}(\extd a))\omega-\cj_{\nabla_{g\uo}g\ut}(\extd a)\omega\\
&=&a\Delta_{LB}\omega+(\Delta_{LB}a)\omega+\nabla_{\extd a+\sigma_{g\uo}(g\ut\tens_A\extd a)}\omega
\end{eqnarray*}
where we used the definition, the bimodule covariant derivative properties and the fact (see the Preliminaries) that the half-curvature $\rho$ depends on the element of $\Omega^1\tens_A\Omega^1$, so that $\rho(g\uo\tens_A g\ut a)=\rho(a g\uo\tens_A g\ut)=a\rho(g)$ by centrality of the metric and the covariant derivative properties. We see that
\begin{equation}\label{LBleib} L_{\Delta_{LB}}(a,\omega)=\nabla_{\extd a+\sigma_{g\uo}(g\ut\tens_A\extd a)}\omega,\quad\forall a\in A,\ \omega\in\Omega\end{equation}
quite generally. Putting in the specific form of $\sigma$ in our case, we immediately obtain the expression stated. Meanwhile,  from (\ref{l0}) we see that $L_\Delta(a,\omega)=2\nabla_{\extd a}\omega$ so of $g^1g^2=0$ then the difference $\Delta_{LB}-\Delta$ is a left $A$-module map and in that case we can define ${\rm Ricci}_{\Delta}$ as stated.  \endproof

This should be viewed as a working definition and novel approach to the Ricci tensor in noncommutative geometry,  motivated by our classical calculations in Section~\ref{secricci}. It does not necessarily connect up to the trace of the Riemann curvature in general, but does do so in the classical Riemannian manifold case.

\subsection{Flat cleft central extensions and their construction}

To proceed further we say that a central extension if {\em flat} if $\Delta$ (but not necessarily the whole cocycle) is cohomologous to zero. This means that up to an isomorphism we can take $\Delta=0$, which is clearly a natural restriction. According to our analysis of morphisms in Lemma~\ref{equivsols}, this is equivalent to the existence of a degree -1 linear map $\delta$ such that $\Delta=\extd\delta+\delta\extd$. Meanwhile, we have seen that a cleft extension is controlled entirely by $\Delta$ and hence now by $\delta$. 

\begin{proposition} In a flat cleft extension, if $\delta$ is regular in the sense of Definition~\ref{regular} then the cleft extension is 1-regular in the sense of Definition~\ref{cleftreg} and $(\ ,\ )$ coincides with the metric associated to $\delta$. \end{proposition}
\proof From (\ref{l0}) we have
\begin{eqnarray*} \lo \extd a,b\lc&=&L_\Delta(a,b)=\Delta(ab)-(\Delta a)b-a\Delta(b)=\delta\extd(ab)-(\delta\extd a)b-a\delta\extd b\\ 
&=&\delta((\extd a)b)+\delta(a\extd b)-(\delta\extd a)b-a\delta\extd b=\extd a\rinter_{\extd b}+\inter_{\extd a}\extd b=2(\extd a,\extd b).\end{eqnarray*}
where $(\ ,\ )$ is from Definition~\ref{regular}. Then from (\ref{l2}) we have 
\[2J_{\extd a}(b\extd c)=\lo(\extd a)b,c\lc=\lo\extd a,bc\lc-\lo\extd a,b\lc c=2(\extd a,\extd(bc))-2(\extd a,\extd b)c=2(\extd a,b\extd c)\]
for all $a,b,c$, using the bimodule properties of $(\ ,\ )$. Hence $\cj$ on degree 1 is well-defined and agrees with $(\ ,\ )$ from $\delta$. 
\endproof

For higher degrees one needs a derivation property as in Lemma~\ref{deriv} which works well in the graded-commutative case covered later, or a comparable assumption in the general context as in the following theorem.

\begin{theorem}\label{construct} Let $\perp$ be a degree -2 bilinear map on  a standard DGA $\Omega(A)$ such that $\perp a=a\perp=0$ for all $a\in A$ and 
\[ (-1)^{|\eta|} (\omega\eta)\perp\zeta+(\omega\perp\eta)\zeta=\omega\perp(\eta\zeta)+(-1)^{|\omega|+|\eta|}\omega(\eta\perp\zeta),\quad\forall\omega,\eta,\zeta\in\Omega\]
and let $\delta$ be regular with
\[ \delta(a\omega)-a\delta\omega=\extd a\perp\omega,\quad \forall a\in A,\ \omega\in \Omega.\]
 Then there is a regular flat cleft extension with 
\[ \Delta=\extd\delta+\delta\extd,\quad \lo\omega,\eta\lc=L_\delta(\omega,\eta)+\omega\perp\extd\eta-(-1)^{|\omega|}\extd\omega\perp\eta-(-1)^{|\omega|}\extd(\omega\perp\eta),\quad\forall \omega,\eta\in \Omega.\]
\end{theorem}
\proof We first observe that special cases of the $\perp$ identity when one of the forms is in degree 0 tell us that $\perp:\Omega\tens_A\Omega\to \Omega$ and that this is a bimodule map. Moreover, 
\[ \Delta_0=0,\quad \lo\ ,\ \lc_0=\omega\perp\extd\eta-(-1)^{|\omega|}\extd\omega\perp\eta-(-1)^{|\omega|}\extd(\omega\perp\eta)\]
provide a flat extension. For this we check 
\begin{eqnarray*} \lo\omega\eta,\zeta\lc&+&\lo\omega,\eta\lc \zeta-\lo\omega,\eta\zeta\lc-(-1)^{|\omega|}\omega\lo\eta,\zeta\lc\\
&=&(\omega\eta)\perp\extd\zeta-(-1)^{|\omega|+|\eta|}\extd((\omega\eta)\perp\zeta)-(-1)^{|\omega|+|\eta|}  (\extd(\omega\eta))\perp\zeta\\
&+&(\omega\perp\extd\eta)\zeta-(-1)^{|\omega|}(\extd(\omega\perp\eta))\zeta-(-1)^{|\omega|}(\extd\omega\perp\eta)\zeta\\
&-&\omega\perp\extd(\eta\zeta)+(-1)^{|\omega|}\extd(\omega\perp(\eta\zeta))+(-1)^{|\omega|}\extd\omega\perp(\eta\zeta)\\
&-&(-1)^{|\omega|}\omega(\eta\perp\extd\zeta)+(-1)^{|\omega|+|\eta|}\omega\extd(\eta\perp\zeta)+(-1)^{|\omega|+|\eta|}\omega(\extd\eta\perp\zeta)\\
&=&(\omega\eta)\perp\extd\zeta-(-1)^{|\omega|+|\eta|}\extd((\omega\eta)\perp\zeta)-(-1)^{|\omega|+|\eta|}((\extd\omega)\eta)\perp\zeta-(-1)^{|\eta|}(\omega\extd\eta)\perp\zeta\\
&+&(\omega\perp\extd\eta)\zeta-(-1)^{|\omega|}\extd((\omega\perp\eta)\zeta)+(-1)^{|\omega|}(\omega\perp\eta)\extd\zeta-(-1)^{|\omega|}(\extd\omega\perp\eta)\zeta\\
&-&\omega\perp((\extd\eta)\zeta)-\omega\perp(\eta\extd\zeta)+(-1)^{|\omega|}\extd(\omega\perp(\eta\zeta))+(-1)^{|\omega|}\extd\omega\perp(\eta\zeta)\\
&-&(-1)^{|\omega|}\omega(\eta\perp\extd\zeta)+(-1)^{|\omega|+|\eta|}\extd(\omega(\eta\perp\zeta))-(-1)^{|\omega|+|\eta|}(\extd\omega)(\eta\perp\zeta)\\
&&\quad\quad+(-1)^{|\omega|+|\eta|}\omega(\extd\eta\perp\zeta)\\
&=&0
\end{eqnarray*}
There are 16 terms and they cancel in groups of 4 under application of the 4-term condition on $\perp$ assumed in the statement of the theorem when applied to appropriate elements. For example, the leading term is in a group of 4 which cancel by application to $\omega,\eta,\extd\zeta$. Next, for any degree -1 map $\delta$ we add its coboundary according to Lemma~\ref{equivsols} to obtain the stated extension $(\Delta,\lo\ ,\ \lc)$. This is cleft iff $\delta$ obeys the condition stated given that $\perp a=a\perp=0$ for all $a\in A$. This is also half of the assumed regularity of $\delta$ (namely that $\inter_{\extd a}=\extd a\perp$). Also, we find $\lo\omega a,b\lc=L_\delta(\omega a,b)+\omega \perp a\extd b)$ so for this to depend only on $a\extd b$ we need the other half of the regularity assumption on $\delta$, namely $L_\delta(\omega a,b)=(\omega a)\rinter_{\extd b}=\omega\rinter_{a\extd b}$. Then we define $ \cj_\omega(a\extd b)={1\over 2}( L_\delta(\omega a,b)+\omega\perp (a\extd b))$, or
\begin{equation}\label{jflatcleft} \cj_\omega(\zeta)={1\over 2}(\omega\rinter_\zeta+\omega\perp\zeta),\quad\forall \omega\in \Omega,\ \zeta\in\Omega^1.\end{equation}
According to Proposition~\ref{cleftcon1}, $\nabla_\omega={1\over 2} \lo \omega,\ \lc$ then gives us a bimodule covariant derivative, in fact extended to a covariant derivative along $\omega$ of all degrees. \endproof

\begin{lemma} \label{perpB} Let $\Omega(A)$ be a standard DGA and $(\perp,\delta)$ a solution for the data in Theorem~\ref{construct}. Then
\[\omega\perp'\eta=\omega\perp\eta+(-1)^{|\omega|+1}L_B(\omega,\eta),\quad\delta'=\delta+B\extd -\extd B\]
is also a solution, for any degree -2 bimodule map $B$. This leaves $\Delta$ and $\lo\ ,\ \lc$ in Theorem~\ref{construct} and hence the induced metric and covariant derivative unchanged. \end{lemma}
\proof This is a matter of direct verification that $\perp'$ still obeys the 4-term relation in Theorem~\ref{construct}. Moreover.
\begin{eqnarray*} \delta'(a\omega)-a\delta'&=&\extd a\perp\omega+B\extd(a\omega)-\extd B(a\omega)-a(B\extd -\extd B)\omega\\
&=&\extd a\perp\omega+B((\extd a)\omega)-(\extd a)B\omega=\extd a\perp'\omega\end{eqnarray*}
for all $a\in A$, $\omega\in\Omega$. Hence the flatness condition is maintained as is the left half of the regularity of $\delta'$. We also have 
\begin{eqnarray*}L_{\delta'}(\omega a,b)&=&L_\delta(\omega a,b)+B\extd(\omega ab)-\extd B(\omega ab)-(B\extd(\omega a))b-(\extd B(\omega a))b\\
&=&L_\delta(\omega a,b)+B(\extd(\omega ab)-(\extd(\omega a))b)-\extd ( (B(\omega a))b))+(\extd B(\omega a))b\\
&=& L_\delta(\omega a,b)+(-1)^{|\omega|}B(\omega a\extd b)-(-1)^{|\omega|}(B \omega)a\extd b
\end{eqnarray*}
using that $B$ is a bimodule map. The last two terms depend only on $a\extd b$ so the other half of the required regularity condition is also maintained with $\rinter'_\omega(\eta)=\rinter_\omega(\eta)+(-1)^{|\omega|}(B(\omega\eta)-B(\omega)\eta)$. Since the extension is cleft, $\lo\ ,\ \lc$ is determined by $\Delta$ and hence is unchanged as the latter is unchanged by the addition of $\extd(B\extd-\extd B)+(B\extd-\extd B)\extd=0$. \endproof

In particular, we can start with the zero solution, then any degree -2 bimodule map $B$ generates a `coboundary' solution but with trivial end product.

\subsection{Noncommutative inner flat cleft extensions}

Before we do the classical case we present a class of `quantum' or noncommutative examples. We focus on the inner case, which is not possible classically, where we assume the existence of a 1-form $\theta\in\Omega^1$ such that $\extd\omega=\theta\omega-(-1)^{|\omega|}\omega\theta$ for all $\omega\in \Omega(A)$, which we assume to be of standard type. 

\begin{proposition}\label{inner} If $\Omega(A)$ is of standard type and inner via $\theta\in\Omega^1$ and if $\perp$ solves the 4-term condition in Theorem~\ref{construct} then $\delta=\theta\perp$ completes the data for a regular flat cleft extension. Here
\[ j_\omega(\zeta)={1\over 2}\omega\perp\zeta,\quad  \Delta =2\nabla_\theta-\theta^2\perp\]
\[  \nabla_\omega=-{1\over 2}L_{\perp\theta}(\omega,\ ),\quad \sigma_\omega(\eta\tens_\zeta)={1\over 2}\left((\omega\eta)\perp\zeta-(-1)^{|\omega|}\omega(\eta\perp\zeta)\right)\]
for all $\omega,\eta\in \Omega$ and $\zeta\in \Omega^1$. The cocycle is $\lo \omega,\ \lc=2\nabla_\omega$. 
\end{proposition} 
\proof We have assumed $\perp$ and clearly $\delta(a\omega)-a\delta\omega=\theta a\perp \omega-a\theta\perp\omega=\extd a\perp\omega$. Similarly $L_\delta(\omega,a)=\delta(\omega a)-(\delta\omega)a=(\theta\perp\omega a)-(\theta\perp\omega)a=0$ by the bimodule properties of $\perp$. Hence $\delta$ is regular and we have a regular flat cleft extension. Clearly $\cj_\omega(\extd a)={1\over 2}\omega\perp\extd a$. To compute the covariant derivative,
\begin{eqnarray*} \omega\perp\extd\eta&-&(-1)^{|\omega|}\extd\omega\perp\eta-(-1)^{|\omega|}\extd(\omega\perp\eta)\\
&=&\omega\perp(\theta\eta)-(-1)^{|\eta|}\omega\perp(\eta\theta)-(-1)^{|\omega|}(\theta\omega)\perp\eta+ (\omega\theta)\perp\eta\\
&&-(-1)^{|\omega|}\theta(\omega\perp\eta)+(-1)^{|\eta|}(\omega\perp\eta)\theta\\
&=&(\omega\perp\theta)\eta-(\omega\eta)\perp\theta-\theta\perp (\omega\eta)+ (-1)^{|\omega|}\omega(\theta\perp\eta)\\
&&+(\theta\perp\omega)\eta+(-1)^{|\omega|}\omega(\eta\perp\theta)\\
&=&-L_{\theta\perp}(\omega,\eta)-L_{\perp\theta}(\omega,\eta)
\end{eqnarray*}
using the definition of $\extd$ in the inner case and applying the 4-term identity to the terms pairwise, 3 times. The generalised braiding is from Proposition~\ref{cleftcon2}. For the Hodge Laplacian, 
\begin{eqnarray*}\Delta\omega&=& \theta(\theta\perp\omega)+(-1)^{|\omega|}(\theta\perp\omega)\theta+\theta\perp(\theta\omega)-(-1)^{|\omega|}\theta\perp(\omega\theta)\\
&=&-\theta^2\perp\omega+(\theta\perp\theta)\omega-(\theta\omega)\perp\theta-\theta(\omega\perp\theta)=-\theta^2\perp\omega-L_{\perp\theta}(\theta,\omega)
\end{eqnarray*}
using the definition of $\Delta$, $\extd$ and two applications of the 4-term identity for $\perp$. \endproof

We have covered the case of $\nabla_\omega$ along forms of all degrees but $\omega\in \Omega^1$ corresponds to a usual covariant derivative.  Note also that $(\ ,\ )={1\over 2}\perp$ in this class of examples and Lemma~\ref{perpB} provides a construction for $\perp$. 

\begin{example}\label{Z2}  We let $A=k(\{x,y\})=k\oplus k$,  the algebra of functions on 2 points and $\Omega(A)$ its universal calculus. Here for any unital algebra, $\Omega^n\subset A^{\tens(n+1)}$ is the sub-bimodule such that the product applied to any two adjacent copies is zero. The exterior derivative is $\extd(a_0\tens\cdots\tens a_n)=\sum_i(-1)^{i}a_0\tens\cdots \tens a_{i-1}\tens 1\tens a_i\tens \cdots\tens a_n$. In practice it in our case is better to think of $A$ as functions on the group $\Z_2$. The universal calculus is on a group is bicovariant and hence has a basis of left-invariant 1-forms, and also the calculus is inner.  In our case it is generated by $A$ and a single 1-form $\theta$, so in degree $n$ the $n$-form $\theta^n$ forms a basis.  If $f\in A$, let $\bar f(x)=f(y)$, $\bar f(y)=f(x)$. Then the relations of the DGA are
 \[ \theta f= \bar f \theta,\quad\extd f=(\bar f-f)\theta,\quad \extd \theta^n= (1-(-1)^{n})\theta^{n+1}.\]

 We next solve the condition in Theorem~\ref{construct}. Since $\perp$ is a bimodule map it is enough to define it on the invariant forms, i.e. on powers of $\theta$, and we take for example 
 \[ \theta^m\perp \theta^n=2 (-1)^{m+1}mn \theta^{m+n-2},\]
where the 2 also fixes a particular normalisation that we will need. 

According to Proposition~\ref{inner} we now have a regular flat cleft extension and the associated codifferential, bimodule covariant derivative, invertible metric and  Laplacian are \[ \delta (f\theta^n)=2\bar f n \theta^{n-1},\quad \nabla_\theta (f \theta^n)=(f-(-1)^n\bar f)\theta^n,\quad \sigma_\theta(f\theta^n\tens f'\theta)=(-1)^n\bar f \theta^n \bar{f'}\]
\[(f\theta,f'\theta)=f\bar{f'},\quad g=\theta\tens \theta,\quad\Delta\omega=2(\nabla_\theta\omega+2 |\omega|\omega),\quad\forall f,f'\in A,\ \omega\in \Omega.\]
These computations are immediate from the general structure in Proposition~\ref{inner} applied in our case. We find now that the connection is torsion free and metric-compatible. Thus, the torsion is
\[ T(f\theta^n)=\theta\nabla_\theta (f\theta^n)-(\bar f-(-1)^nf)\theta^{n+1}=0,\]
while from the action (\ref{covtens}) of tensor products we have
\[ \nabla_\theta g=\nabla_\theta(\theta\tens\theta)=\nabla_\theta\theta\tens\theta+\sigma_\theta(\theta\tens\theta)\tens_A\nabla_\theta\theta=0\]
using $\nabla_\theta=2\theta$ and $\sigma_\theta(\theta\tens\theta)=-\theta$.

 One also sees that $\delta^2$ is not zero but commutes with functions and that, more surprisingly, the canonical Laplace-Beltrami operator vanishes,
 \[ \Delta_{LB}=\nabla_{g^1}\nabla_{g^2}-\nabla_{\nabla_{g^1}g^2}=(\nabla_\theta)^2- 2\nabla_\theta=0.\]
 This happens to be tensorial in our example, which is possible as the metric is not quantum symmetric. One also finds that the noncommutative Riemann curvature vanishes. Using the
 algebraic form, this is 
 \[ R(\theta^n)=(\extd\tens\id-(\wedge\tens\id)(\id\tens\nabla))\nabla \theta^n=\extd 2 \theta\tens\theta^n-2\theta2\theta\tens\theta^n=0\]
 for all $n$ odd (and zero in any case if $n$ is even).  
 
 Finally, one may similarly compute from $L_{\perp\theta}(\theta^m,\theta^n)$ that
 \[ \nabla_{\theta^m}(\theta^n)= (-1)^{m+1}m\theta^{m-1}\nabla_\theta=j_{\theta^m}(\theta)\nabla_\theta\]
 so that the braided-Leibniz rule applies by Proposition~\ref{braleib}. One  can also verify this  directly as a useful check. 
 \end{example}

\subsection{Classical type flat cleft extensions}

Finally we specialise to the case where $\Omega(A)$ is graded-commutative case and of `classical type'. First we consider the general theory of cleft extensions of classical type.

\begin{proposition}\label{cleftclass} For a regular cleft extension, on $\Omega(A)$ of classical type, the associated $(\ ,\ )$ is symmetric, the interior product  is a graded derivation and the covariant derivative in Proposition~\ref{cleftcon1} has symmetric part
\[ \nabla_\omega\eta+\nabla_\eta\omega=\cj_{\extd\omega}(\eta)+\cj_{\extd \eta}(\omega)+\extd(\omega,\eta),\quad\forall\omega,\eta\in \Omega^1\]
and has torsion and metric-compatibility tensor obeying
\[T(\zeta,\omega)(\eta)+T(\zeta,\eta)(\omega)=C_\zeta(\omega,\eta),\quad\forall \omega,\eta,\zeta\in \Omega^1.\]
In the invertible case the Laplace-Beltrami operator obeys 
\[ L_{\Delta_{LB}}(a,\omega)=2\nabla_{\extd a}\omega,\quad a\in A,\ \omega\in \Omega\]
so that ${\rm Ricci}_\Delta$ is defined. 
\end{proposition}
\proof From (\ref{l0}) we have $L_\Delta(a,b)=\lo \extd a,b\lc=2(\extd a,\extd b)$ which is symmetric by graded-commutativity. Next, comparing (\ref{l2}),(\ref{l3}) and using the graded-commutativity and (\ref{l1}) we see that $\lo\ ,a\lc$ and hence $\cj_{(\ )}(\extd a)$ is a graded derivation. This is also the map $\inter_{\extd a}$ in the general theory of form-covariant derivatives in Section~2.2. 
Similarly, from (\ref{l00}) and (\ref{l0}) and  graded-commutativity we have
\[ \nabla_{\extd a}\omega={1\over 2}L_\Delta(a,\omega)={1\over 2}L_\Delta(\omega,a)=-\nabla_\omega\extd a+\cj_{\extd\omega}(\extd a)+\extd \cj_\omega(\extd a)\]
from which we conclude the stated symmetry of the covariant derivative. The torsion result is then immediate from Lemma~\ref{symanti} after allowing for the change of notation.  Also, since $g$ is symmetric, the Laplace-Beltrami operator in Proposition~\ref{cleftricci} has the stated Leibnizator (same as the Hodge Laplacian). \endproof

The generalised covariant derivatives along higher forms also apply as in the general case, but with $\sigma$ trivial in the sense $\sigma_\omega(\eta\tens_A\zeta)=\cj_\omega(\zeta)\eta$. Thus, 
\begin{corollary}\label{higher} For  a regular cleft extension on $\Omega(A)$ of classical type the `extended covariant derivative' on $\Omega$,
\[ \nabla_\omega:={1\over 2}\lo\omega,\ \lc\]
has degree $|\omega|-1$, obeys 
\[ \nabla_{a\omega}=a\nabla_\omega,\quad \nabla_\omega(a\eta)=a\nabla_\omega \eta+\cj_\omega(\extd a)\eta,\quad\nabla_1\eta=\nabla_\omega 1=0,\]
\[ {1\over 2}L_\Delta(\omega,\ )=[\extd,\nabla_\omega\}+\nabla_{\extd\omega},\quad\forall a\in A,\ \omega,\eta\in \Omega,\]
and if $(\ ,\ )$ is invertible and if $\nabla_\omega$ is a derivation for all $\omega\in \Omega^1$ then $\nabla_\omega$ is a graded-derivation for all $\omega\in\Omega$ and 
\[ \nabla_{\omega_1\cdots\omega_m}=\sum_{i=1}^m(-1)^{i-1}\omega_1\cdots\widehat{\omega_i}\cdots\omega_m\nabla_{\omega_i},\quad\forall \omega_i\in
\Omega^1,\]
where the hat denotes omission.
\end{corollary}
\proof We bring together some of the properties of the solution for $\lo\ ,\ \lc$ obtained in the course of Theorem~\ref{cleftclasslevi}. These include (\ref{c1})-(\ref{c2}) by definition of  $\nabla_\omega={1\over 2}\lo\omega,\ \lc$. 
 In this case the explicit formula follows from 
\[ \nabla_{\omega\eta}=(-1)^{(|\omega|-1)|\eta|}\eta\nabla_\omega+(-1)^{|\omega|}\omega\nabla_\eta,\quad \forall \omega,\eta\in \Omega\]
deduced from (\ref{c1}) when $\nabla_\omega$ is a  graded-derivation. We use Proposition~\ref{braleib}.    \endproof

We now consider the converse direction. In Section~2 we had a notion of $\delta$ of classical type and we assume this now except {\em without} the $\delta^2$ tensorial and {\em without} the
symmetry of $\delta$. Thus we assume that $\delta$ is regular in the sense of Definition~\ref{regular} and that the conditions in part (2) of Lemma~\ref{deriv} apply.

\begin{theorem}\label{cleftclasslevi}  Let $\Omega(A)$ be of classical type and $\delta$ a regular degree -1 map obeying the derivation conditions in Lemma~\ref{deriv}. Then $\perp$ defined on degree 1 by $\extd a\perp\omega=\delta(a\omega)-a\delta\omega$  extends as a graded-derivation of appropriate degree and Theorem~\ref{construct} provides a regular flat cleft extension
\[ \Delta=\extd\delta+\delta\extd,\quad  \lo\omega,\eta\lc=L_\delta(\omega,\eta)+\CL_\omega\eta- (-1)^{|\omega|}(\extd\omega)\perp\eta,\quad\forall \omega,\eta\in \Omega(A)\]
with associated metric $(\omega,\eta )={1\over 2}(\omega\perp\eta+\eta\perp\omega)$ and $\nabla_\omega={1\over 2}\lo\omega,\ \lc$. Here
\[ \CL_\omega\eta=\omega\perp\extd\eta-(-1)^{|\omega|}\extd(\omega\perp\eta),\quad\forall \omega,\eta\in \Omega\]
extends the usual Lie derivative as a degree $|\omega|-1$ derivation. The above provides the unique cleft central extension with the given $\Delta$. If $\delta$ is of classical type then the covariant derivative is torsion free (and hence metric compatible) and is a derivation. 
\end{theorem}
\proof We let $\perp:\Omega^1\tens_A\Omega^1\to A$ be defined by $\delta(a\omega)-a\delta\omega=\extd a\perp\omega$ (so that $\eta\perp\omega=\inter_\eta(\omega)$ for $\omega,\eta\in \Omega^1$ in the notation of Section~2.1).  This extends to 
to a map 
\[ \omega_1\cdots\omega_m\perp \eta_1\cdots\eta_n=\sum_{i,j}(-1)^{i+j}(\omega_i\perp\eta_j)\omega_1\cdots\widehat{\omega_i}\cdots\omega_m\eta_1\cdots\widehat{\eta_j}\cdots\eta_n\]
much as in Section~2. This is antisymmetric in the $\omega_i$ factors and the $\eta_i$ factors, and hence is well-defined on $\Omega\tens_A\Omega$. Also as in Section~2,  $\omega\perp$ is a degree $|\omega|-2$ derivation and $\perp\eta$ similarly obeys
\begin{equation}\label{perpeta} (\omega\omega')\perp\eta=(-1)^{|\omega'|(|\eta|-1)}(\omega\perp\eta)\omega'+(-1)^{|\omega|}\omega(\omega'\perp\eta),\quad\forall\omega,\omega',\eta\in \Omega. \end{equation}
Then the 4-term relation required in  Theorem~\ref{construct} holds as
\begin{eqnarray*} &&\kern-10pt (-1)^{|\eta|}(\omega\eta)\perp\zeta-(-1)^{|\omega|+|\eta|}\omega(\eta\perp\zeta)=(-1)^{|\eta|}(-1)^{|\eta|(|\zeta|-1)}(\omega\perp\zeta)\eta\\
&&=(-1)^{|\eta||\zeta|}(\omega\perp\zeta)\eta=(-1)^{|\eta||\omega|}(-1)^{|\eta|(|\omega|+|\zeta|-2)}(\omega\perp\zeta)\eta\\
&&=(-1)^{|\eta|(|\omega|-2)}\eta(\omega\perp\zeta)=\omega(\eta\perp\zeta)-(\omega\perp\eta)\zeta.\end{eqnarray*}
Next, by the derivation assumption in Lemma~\ref{deriv} we know that $\delta$ is such that $L_\delta(a,\omega)=\extd a\perp\omega$ for all $\omega\in\Omega$, as needed for flatness. This also meets the regularity requirement and we have a regular flat cleft extension. The interior product is $\cj_\omega(\extd a)={1\over 2}(L_\delta(\omega,a)+\omega\perp\extd a)={1\over 2}(L_\delta(a,\omega)+\omega\perp\extd a)$ due to the graded-commutativity. Thus 
\begin{equation}\label{cjclass}\cj_\omega(\zeta)={1\over 2}(\zeta\perp\omega+\omega\perp\zeta),\quad \omega\in\Omega,\ \zeta\in\Omega^1.\end{equation}
Using (\ref{perpeta}), it follows that $\cj_{(\ )}(\zeta)$ is a degree -1 derivation for all $\zeta\in \Omega^1$ extending the metric $(\ ,\ )$. It therefore concides with the map $\inter_\zeta$ in Section 2.2 which developed the general theory of covariant derivatives on $\Omega(A)$ of classical type.  The  covariant derivative also has precisely the form in Theorem~\ref{levi} but generalises it as we have not required the symmetry condition on $\delta$ nor that  $\delta^2$ to be tensorial. Because we have not assumed the symmetry condition, the map $\inter$ associated to $\delta$ in Section~2.1 is no longer the interior product associated to the metric as used in Section~2.2; we are denoting the latter now as $\cj$.  When we do have $\delta$ of classical type we see that $(\ ,\ )=\perp$ on degree 1 and the theory of Section~2 applies directly, with $\perp$ now having the same meaning as in Section~2.1. Uniqueness is from Lemma~\ref{uniquesol}.

Finally, we prove that 
\[ \CL_\omega(\eta\eta')=(\CL_\omega\eta)\eta'+(-1)^{(|\omega|-1)|\eta|}\eta\CL_\omega\eta',\quad\forall \omega,\eta,\eta'\in \Omega.\]
This follows straightforwardly from 
\[ \CL_\omega(\eta\eta')=\omega\perp\left((\extd\eta)\eta'+(-1)^{|\eta|}\eta\extd\eta'\right)-(-1)^{|\omega|}\extd(\omega\perp(\eta\eta'))\]
on computing further via  $\omega\perp$ a degree $|\omega|-2$ derivation and then comparing with the right hand side computed from the definition. This perhaps justifies our term `extended Lie derivative'. Note that $\CL_\omega=\omega\perp\extd - (-1)^{|\omega|}\extd \omega\perp$  reduces to the Cartan formula for the usual Lie derivative when $\omega\in\Omega^1$.  \endproof

This puts the Levi-Civita covariant derivative of Section~2 into a slightly wider context coming from the theory of flat cleft extensions. This more general theory applies to more general BV structures on $\Omega(A)$ of classical type as we do not assume that $\delta$ is symmetric. Note that by using the freedom in Lemma~\ref{perpB} one could change the antisymmetric part of $\perp$ on degree 1 and hence potentially make it symmetric but not necessarily obtaining or retaining other desired properties of $\delta$. 

The uniqueness in Theorem~\ref{cleftclasslevi} is an analogue of the way that in Riemannian geometry the Levi-Civita covariant derivative is uniquely determined by the metric, in the present case encoded in the choice of $\perp$ in degree 1 or more precisely by $\delta$ giving this. Theorem~\ref{cleftclasslevi} also achieves our goal of putting in a proper context formulae in \cite[Sec.~2]{Ma:bh} where $A=C(M)$ the smooth functions (say) on a Riemannian manifold  $(M,g)$. 

\begin{corollary}\label{Riequant} Every classical Riemannian manifold $M$ has an `almost commutative' exterior algebra $\tilde{\Omega}(M)$ with extension data
\[ \lo a,\eta\lc =0,\quad  \lo \omega,a \lc=2(\extd a,\omega),\quad \lo \omega,\eta\lc=2\nabla_\omega\eta,\quad\forall a\in C(M),\ \omega,\eta\in \Omega^1(M)\]
and  $\Delta$ the Hodge Laplacian. We obtain the relations and differential
\[ [a,\omega]=\lambda(\extd a,\omega)\theta',\quad\{\omega,\eta\}=\lambda(\CL_{\omega}\eta+\inter_\eta\extd\omega)\theta'\]
\[  [a,\theta']=0,\quad \{\omega,\theta'\}=0,\quad \theta'^2=0\]
\[ \extd_\cdot a= \extd a- {\lambda\over 2} (\Delta_{LB} a)\theta',\quad\extd_\cdot \omega=\extd\omega+{\lambda\over 2}((\Delta_{LB}-{\rm Ricci})\omega)\theta'\]
for all $a\in C(M),\ \omega,\eta\in\Omega^1$, essentially as in  \cite[Sec.~2]{Ma:bh} in the case $\extd\theta'=0$.
\end{corollary}
\proof This follows from Theorem~\ref{cleftclasslevi} specialised to the classical case as in Remark~\ref{classical}. The $\cdot$ product and $\extd_\cdot$ from Proposition~\ref{Omegaext} then give the commutation relations and derivative as stated. This is not quite the generality in \cite{Ma:bh} where we did not assume that $\extd\theta'=0$ (we come to this later) and we have used (\ref{connpm}) to simplify in terms of the Lie derivative. There is also a change of sign of $\lambda$ compared to \cite{Ma:bh} and we did not give $\extd_\cdot\omega$ so explicitly as we do now. The conversion from the Hodge-Laplacian to the Laplace-Beltrami is standard but note that by the 2nd order Leibniz rule in Remark~\ref{classical} the leading part of $\Delta$ in our conventions agrees with the Laplace-Beltrami operator, after which the coefficient of Ricci can be fixed by the identity $[\Delta_{LB},\extd]a={\rm Ricci}(\extd a)$ used in \cite{Ma:bh} as being equivalent to $[\Delta,\extd]=0$.   \endproof 

The work \cite{Ma:bh} then `quantises' the Levi-Civita covariant derivative to a bimodule covariant derivative on this extended DGA. This comes with a generalised braiding $\sigma$ for which the failure to obey the braid relations expresses the classical Riemann curvature and for which the failure to obey $\sigma^2=\id$ expresses the classical Ricci curvature. We refer to \cite[Sec 2]{Ma:bh} for more details. 
 
\begin{remark}  In the proof of Theorem~\ref{construct} we solved first for a flat but not cleft central extension which would now be
\begin{equation}\label{Omegadelta2}\Delta=0,\quad \lo\omega ,\eta \lc=\CL_\omega\eta- (-1)^{|\omega|}(\extd\omega)\perp\eta,\quad\forall \omega,\eta\in \Omega(A).\end{equation}
Thus, starting with a classical manifold we obtain an isomorphic quantisation to Corollary~\ref{Riequant}, but with the merit of not deforming $\extd$. We view this as a base noncommutative geometry of which Corollary~\ref{Riequant} is then a `stochastic differential' version.  It commutation relations are necessarily the same as in Corollary~\ref{Riequant}. \end{remark}

\subsection{Further extension of a flat cleft extension}

Here to cover the case $\extd\theta'\ne 0$, which turns out to be necessary for later sections. 

\begin{proposition}\label{Omegafull} Let $\Omega(A)$ be a standard DGA. Its flat central extension constructed in Theorem~\ref{construct} has a further extension $\tilde{\tilde{\Omega}}(A)=\Omega(A)\oplus\Omega(A)\theta'\oplus\Omega(A)\extd\theta'$ with product and exterior derivative
\[ \omega\cdot\eta=\omega\eta+{\lambda\over 2}(-1)^{|\omega|+|\eta|}\lo\omega,\eta\lc\theta'-{\lambda\over 2}(-1)^{|\omega|}(\omega\perp\eta)\extd\theta' \]
\[ \extd_\cdot\omega=\extd \omega - {\lambda\over 2}(-1)^{|\omega|}(\Delta\omega)\theta'+{\lambda\over 2}(\delta\omega)\extd\theta',\quad\forall \omega,\eta\in \Omega(A),\]
where $\theta'^2=\theta'\extd\theta'=(\extd\theta')\theta'=\{\omega,\theta'\}=0$. Moreover $\tilde{\tilde{\Omega}}(A)\to \tilde{\Omega}(A)\to \Omega(A)$ are surjections of DGAs successively setting $\extd\theta'=0$ and $\theta'=0$. 
\end{proposition}
\proof We define of course $\extd_\cdot\theta'=\extd\theta'$ as the notation suggests. For the Leibniz rule we recompute the proof of Proposition~\ref{Omegaext}, 
\begin{eqnarray*}\extd_\cdot(\omega\cdot\eta)&=&\extd_\cdot\left( \omega\eta+{\lambda\over 2}(-1)^{|\omega|+|\eta|}\lo\omega,\eta\lc\theta'-{\lambda\over 2}(-1)^{|\omega|}(\omega\perp\eta)\extd\theta' \right) \\
&=&\extd(\omega\eta)-{\lambda\over 2}(-1)^{|\omega|+|\eta|}(\Delta(\omega\eta))\theta'+{\lambda\over 2}(\delta(\omega\eta))\extd\theta'\\
&&+{\lambda\over 2}(-1)^{|\omega|+|\eta|}(\extd \lo\omega,\eta\lc\theta'+(-1)^{|\omega|+|\eta|-1}\lo\omega,\eta\lc\extd\theta')\\
&&-{\lambda\over 2}(-1)^{|\omega|}(\extd(\omega\perp\eta))\extd\theta' 
\end{eqnarray*}\begin{eqnarray*}
(\extd_\cdot\omega)\cdot\eta&+&(-1)^{|\omega|}\omega\cdot\extd_\cdot\eta=\left(\extd\omega-{\lambda\over 2}(-1)^{|\omega|}(\Delta\omega)\theta'+{\lambda\over 2}(\delta\omega)\extd\theta'\right)\cdot\eta\\
&&\quad\qquad\qquad+(-1)^{|\omega|}\omega\cdot\left(\extd \eta-{\lambda\over 2}(-1)^{|\eta|}(\Delta\eta)\theta'+{\lambda\over 2}(\delta\eta)\extd\theta'\right)\\
&=&(\extd\omega)\eta+(-1)^{|\omega|}\omega\extd\eta+{\lambda\over 2}(-1)^{|\omega|+1+|\eta|}\lo\extd\omega,\eta\lc\theta'+{\lambda\over 2}(-1)^{|\eta|+1}\lo\omega,\extd\eta\lc\theta'\\
&&-{\lambda\over 2}(-1)^{|\omega|+1}(\extd\omega)\perp\eta\extd\theta'-{\lambda\over 2}(\omega\perp\extd\eta)\extd\theta'\\
&&-{\lambda\over 2}(-1)^{|\omega|+|\eta|}(\Delta\omega)\eta\theta'- {\lambda\over 2}(-1)^{|\eta|+|\omega|}\omega(\Delta\eta)\theta'\\
&&+{\lambda\over 2}(\delta\omega)\eta\extd\theta'+{\lambda\over 2}(-1)^{|\omega|}\omega\delta\eta\extd\theta'
\end{eqnarray*}
where we used $\theta'^2=0$ and $\theta'\extd\theta'=0$. We already have equality for all terms except those with $\extd\theta'$ which had been ignored before. Equality of these is exactly the formula for $\lo\omega,\eta\lc$ in Theorem~\ref{construct}. We also need 
\begin{eqnarray*}\extd_\cdot\extd_\cdot\omega&=&\extd_\cdot(\extd \omega-{\lambda\over 2}(-1)^{|\omega|}(\Delta\omega)\theta'+{\lambda\over 2}(\delta\omega)\extd\theta')\\
&=& \extd^2\omega-{\lambda\over 2}(-1)^{|\omega|+1}(\Delta\extd\omega)\theta'+{\lambda\over 2}(\delta\extd\omega)\extd\theta'\\
&&-{\lambda\over 2}(-1)^{|\omega|}((\extd\Delta\omega)\theta'+(-1)^{|\omega|}(\Delta\omega)\extd\theta')+{\lambda\over 2}(\extd\delta\omega)\extd\theta'=0\end{eqnarray*}
where the new $\extd\theta'$ terms cancel in view of the definition of $\Delta$ as the Hodge Laplacian. 
Last but not least, we have to check associativity for the new product. Referring to the proof of Proposition~\ref{Omegaext} the new terms are 
\[ (\omega\cdot\eta)\cdot\zeta=\cdots -{\lambda\over 2}(-1)^{|\omega|+|\eta|}(\omega\eta)\perp\zeta\extd\theta'-{\lambda\over 2}(-1)^{|\omega|}(\omega\perp\eta)\zeta\extd\theta'\]
\[ \omega\cdot(\eta\cdot\zeta)=\cdots-{\lambda\over 2}(-1)^{|\omega|}\omega\perp(\eta\zeta)\extd\theta'-{\lambda\over 2}(-1)^{|\eta|}\omega(\eta\perp\zeta)\extd\theta'\]
Equality holds by exactly the 4-term identity for $\perp$ in Theorem~\ref{construct}.  \endproof

When specialised to $\Omega(M)$ in the case of a smooth manifold as in Corollary~\ref{Riequant}, we have a DGA $\tilde{\tilde{\Omega}}(M)$ with relations that recover the general case of \cite[Sec. 2]{Ma:bh} at low degree and which are now given for all degrees.

\section{Semidirect products of differential graded algebras}

In this section we want to apply some of the theory above. Most of the work will be in classical type, so basically classical Riemannian geometry but done using our codifferential approach. We start with some generalises that include the non-graded-commutative case.

\subsection{Semidirect product by the DGA in one variable} 

Let $(\Omega(A),\extd)$ be a DGA equipped with a derivation $\tau$, i.e. a derivation of $\Omega(A)$ as an
algebra which respects degree and commutes with $\extd$. Let $\Omega(t,\extd t)$ be the general bicovariant calculus of the additive line. This is a super-Hopf algebra with
a parameter $\lambda$ and relations and coproduct
\[ [\extd t,t]=\lambda\extd t,\quad \extd t\wedge\extd t=0,\quad \Delta t=t\tens 1+1\tens t,\quad \Delta \extd t=\extd t\tens 1+1\tens\extd t.\] 
We denote by $A_t=A\lcross k[t]$ the semidirect product of $A$ by $t$ with relations $[t,a]=\lambda \tau(a)$ for all $a\in A$. We recall that a DGA is inner if there is a distinguished 1-form $\theta$ that generates $\extd$ by graded-commutator. 

\begin{proposition}\label{semicalc} Given a derivation $\tau$ commuting with $\extd$,  the super-Hopf algebra $\Omega(t,\extd t)$ acts on $\Omega(A)$ and the super-semidirect product $\Omega(A_t)=\Omega(A)\lcross\Omega(t,\extd t)$ is an inner DGA. The relations and exterior derivative are
\[ [t,\omega]=\lambda (\tau-|\omega|)\omega,\quad\forall \omega\in \Omega(A),\quad [\extd t,\  \}=\lambda\extd \]
where we use the graded-commutator.
\end{proposition}
\proof We define an action of $t$ and $\extd t$ by 
\[ t\la \omega=\lambda(\tau - |\omega |)\omega,\quad \extd t \la \omega= \lambda\extd \omega,\quad\forall \omega\in \Omega(A).\]
Clearly $t$ acts as a derivation and $\extd t$ as a graded-derivation (since $\extd$ does). Moreover  $[\extd t, t]=\extd t\la \lambda (\tau(\omega)-|\omega|\omega)-t\la\lambda\extd \omega=\lambda^2(\extd \tau(\omega)-|\omega|\extd\omega)-\lambda^2(\tau(\extd \omega)- (|\omega|+1)\extd\omega)=\lambda^2\extd\omega=\lambda\extd t\la\omega$ hence we have an action of $\Omega(t,\extd t)$ on $\Omega(A)$ as a super-module algebra. Hence the semidirect product is a super-algebra. Its algebra structure is 
\[ (\omega\phi)(\eta\psi)=(-1)^{|\phi\t| |\eta|}\omega(\phi\o\la\eta)\phi\t\psi,\quad\forall\omega,\eta\in\Omega(A),\ \phi,\psi\in \Omega(t,\extd t).\]
where $\phi\o\tens\phi\t$ is the super-coproduct of $\Omega(t,\extd t)$. Putting in the form of the coproduct, we have the relations as stated. More explicitly, writing $\Omega(A_t)=\Omega(A)k[t]\oplus\Omega(A)k[t]\extd t$ and keeping everything normal ordered with $t,\extd t$ to the right, we have
\[ \omega\phi(t)\eta\psi=\omega(\eta\phi(\lambda(\tau-D)+t))\psi\]
\[ \omega(\phi(t)\extd t)\eta\psi=(-1)^{|\eta|}\omega(\eta\phi(\lambda(\tau-D)+t))(\extd t)\psi+\lambda\omega((\extd\eta)\phi(\lambda(\tau-D)+t))\psi\] where the $\lambda(\tau-D)$ is understood to act to the {\em left} and $D$ is the degree operator.  Finally, we define $\extd$ as the graded-commutator in this algebra with $\theta=\extd t$ and verify that this reduces to $\extd t$ on $t$ and $\extd\omega$ on $\omega\in \Omega(A)$, i.e. extends the given ones. Equivalently, we specify
\[ \extd(\omega \phi)=(\extd\omega)\phi + (-1)^{|\omega|}\omega\extd\phi,\quad \omega\in \Omega(A),\ \phi\in \Omega(t,\extd t).\]
One may verify that this defines a graded-derivation. 
\endproof

\subsection{Conformal 1-forms} 

Here we study the notion of a `conformal Killing' vector field or rather 1-form, $\tau$, on $(\Omega,\delta)$ of classical type as in Section~2. In terms of the metric inner product the usual notion that the metric is invariant up to scale is the `metric-conformal' condition
\begin{equation}\label{confmet} (\CL_\tau\omega,\eta)+(\omega,\CL_\tau\eta)=(\tau,\extd(\omega,\eta))+\alpha(\omega,\eta),\quad\forall \omega,\eta\in \Omega^1,\end{equation}
for some function $\alpha\in A$. Classically this is equivalent to the `conformal Killing' condition, e.g.  \cite{Yano}, 
\begin{equation}\label{confkilling}  \nabla_\omega\tau={1\over 2}\inter_\omega\extd\tau+{\alpha\over 2}\omega,\quad\forall \omega\in\Omega^1.\end{equation}
in terms of the Levi-Civita covariant derivative, which is useful for extending the concept to higher degree forms. In both cases the value of $\alpha$  classically is determined. The special case with $\alpha=0$ corresponds in classical geometry to a Killing 1-form. In our case we are taking a coderivation $\delta$ as the starting point and $\nabla$ from Theorem~\ref{levi}. At this level we define:
\begin{definition} \label{liedelta} Let $(\Omega(A),\delta)$ be of classical type and let $D$ be the degree operator. We say that a 1-form  $\tau\in \Omega^1$ is $\delta$-{\em conformal}  if there is some $\alpha\in A$ and some constant $\beta$ such that
\[ [\delta,\CL_\tau]=\alpha\delta + \inter_{\extd\alpha}(D-\beta)\]
holds acting on degree $D=1$. We say that $\tau$ is {\em strongly $\delta$-conformal} if this holds acting on all of $\Omega$.
\end{definition}

The motivation and relationship between these concepts is clarified by the following lemma.

\begin{lemma}\label{metrictau} Let $(\Omega(A),\delta)$ be of classical type and $\tau\in\Omega^1$.

(1)  $\tau$ conformal Killing  for some function $\alpha$ implies $\tau$ metric-conformal and if $(\ ,\ )$ is nondegenerate then these two notions are equivalent.

(2) $\tau$ conformal Killing implies $\CL_\tau=\alpha\tau$ and if $(\ ,\ )$ is invertible then $\alpha=\delta_\nabla\tau/\beta$ and $\beta={1\over 2}(\ ,\ )(g)$.

(3)  $\tau$ conformal Killing with $\alpha=\delta\tau/\beta$ for some constant $\beta$ and  $\delta^2=0$ implies $\tau$ is $\delta$-conformal. 

(4) $\tau$ $\delta$-conformal for some function $\alpha$ and some constant $\beta$ implies $\tau$ is metric-conformal with factor $\alpha$, and in the nondegenerate case, $(\tau,\extd(\alpha-\delta\tau/\beta))=0$. 
\end{lemma}
\proof   (1) From our formula for the Levi-Civita covariant derivative in Theorem~\ref{levi} the requirement of a conformal Killing form in the sense (\ref{confkilling}) becomes
\begin{equation}\label{confkilling1} \delta(\omega\tau)-(\delta\omega)\tau+\CL_\tau\omega=\omega (\alpha-\delta\tau),\quad\forall\omega\in\Omega^1.\end{equation}
From the formula for $\inter_\zeta L_\delta$ in Lemma~\ref{Schouten} we have 
\[ (\zeta,\omega)\alpha=\inter_\zeta L_\delta(\omega,\tau)+(\zeta,\CL_\tau\omega)=(\omega, \extd\inter_\tau\zeta+ \inter_\tau\extd\zeta)-\inter_\tau\extd\inter_\omega\zeta+(\zeta,\CL_\tau\omega)\]
for all $\omega,\zeta\in \Omega^1$, which is (\ref{confmet}). If $(\ ,\ )$ is nondegenerate then we can push this backwards to conclude (\ref{confkilling1}). 

(2) Setting $\omega=\tau$ in (\ref{confkilling1}), we have $\CL_\tau\tau=\alpha\tau$. If $(\ ,\ )$ is invertible then using (\ref{confkilling}) we have 
\[\delta_\nabla(\tau)=\inter_{g\uo}\nabla_{g\ut}\tau={1\over 2}\left(\inter_{g\uo}\inter_{g\ut}\extd\tau+{\alpha}\inter_{g\uo}g\ut\right).\]
The first term vanishes by symmetry of $g$, giving the value of $\alpha$. This assumes $\Omega(A)$ is of classical type but does not refer
to our  given $\delta$. 

(3) Assuming the stated form of $\alpha$,  the equation (\ref{confkilling2}) becomes
\begin{equation*}\label{confkilling2} \delta(\omega\tau)-(\delta\omega)\tau+\CL_\tau\omega=\omega \alpha(1-\beta),\quad\forall\omega\in\Omega^1\end{equation*}
and we apply $\delta$ to this to give  
\[ -(\delta\omega)\delta\tau-\inter_{\tau}\extd(\delta\omega)+\delta\CL_\tau\omega=\alpha(1-\beta)\delta\omega +(1-\beta)\inter_{\extd\alpha}\omega\]
for all $\omega\in \Omega^1$, which is the condition to be $\delta$-conformal provided $\delta\tau=\alpha\beta$.

(4) Now using the condition to be $\delta$-conformal in Definition~~\ref{liedelta}, we compute for any $a\in A$, 
\begin{eqnarray*} &&\delta\CL_\tau(a\omega)-\CL_\tau\delta(a\omega)-\alpha\delta(a\omega)-(1-\beta)\inter_{\extd \alpha}(a\omega)\\
&=&\delta((\tau,\extd a)\omega+a\CL_\tau\omega)-\CL_\tau(a\delta\omega+\inter_{\extd a}\omega)-a\alpha\delta(\omega)-\alpha\inter_{\extd a}\omega-a(1-\beta)\inter_{\extd\alpha}\omega\\
&=&(\tau,\extd a)\delta\omega+\inter_{\extd(\tau,\extd a)}\omega+\inter_{\extd a}\CL_\tau\omega-(\tau,\extd a)\delta\omega-\CL_\tau\inter_{\extd a}\omega-\alpha\inter_{\extd a}\omega\\
&=&\inter_{\extd(\tau,\extd a)}\omega+\inter_{\extd a}\CL_\tau\omega -\CL_\tau\inter_{\extd a}\omega-\alpha\inter_{\extd a}\omega
\end{eqnarray*}
where we used Lemma~\ref{deriv}. Now the initial expression $=0$ if $\tau$ is $\delta$-conformal as $a\omega$ is another 1-form, so we deduce that 
\[ [\inter_{\extd a},\CL_\tau]\omega=\alpha\inter_{\extd a}\omega-\inter_{\CL_\tau\extd a}\omega,\quad \forall a\in A,\ \omega\in\Omega^1.\]
We then consider this same identity for $\eta=\sum b_i\extd a_i$. Using that $\inter$ is $A$-linear and the Leibniz property of $\CL_\tau$ we deduce that  the identity holds for all $\eta\in\Omega^1$ in place of $\extd a$, which is (\ref{confmet}). In the nondegenerate case we conclude by (1) and (2) that $\CL_\tau\tau=\alpha\tau$ and in this case  the condition for a $\delta$-conformal 1-form  and Lemma~\ref{deriv} tells us that 
\[ \inter_{\extd\alpha}\tau=\delta(\alpha\tau)-\alpha\delta\tau=\delta\CL_\tau\tau-\alpha\delta\tau=(\tau,\extd \delta\tau)+(1-\beta)\inter_{\extd\alpha}\tau\]
or $(\tau,\extd(\beta\alpha))=(\tau,\extd\delta\tau)$ which we write suggestively as stated. \endproof

It follows that in the case of a classical Riemannian manifold all three notions are equivalent. This is well-known for (\ref{confmet})-(\ref{confkilling}) whereas our notion of $\delta$-conformal in Definition~\ref{liedelta} is the one we will want and appears to be new. We will need several more properties of $\delta$-conformal 1-forms as follows. 

\begin{lemma}\label{lieperp} Let $(\Omega(A),\delta)$ be of classical type and $\tau$ $\delta$-conformal.  Then 
\begin{equation*} \CL_\tau(\omega\perp\eta)+\alpha\omega\perp \eta=(\CL_\tau\omega)\perp\eta+\omega\perp\CL_\tau\eta\end{equation*}
\begin{eqnarray*} \CL_\tau S(\omega,\eta)+\alpha S(\omega,\eta)-S(\CL_\tau\omega,\eta)&-&S(\omega,\CL_\tau\eta)\nonumber=(-1)^{|\omega|}(\extd\alpha)(\omega\perp\eta)\end{eqnarray*}
for all $\omega,\eta\in \Omega$. Here $S(\omega,\eta)=\omega\perp\extd \eta-(-1)^{|\omega|}\extd(\omega\perp\eta)-(-1)^{|\omega|}(\extd\omega)\perp\eta$.
\end{lemma}
\proof (1) The smallest nontrivial case of the first statement is with $\omega,\eta$ of degree 1 and is just the metric-conformality (\ref{confmet}), so this case holds by Lemma~\ref{metrictau}. We let $\eta\in \Omega^1$ so that $\perp\eta=\inter_\eta$, then we want to prove 
\begin{equation}\label{lieint}
[\inter_\eta,\CL_\tau]=\alpha\inter_\eta-\inter_{\CL_\tau\eta}\end{equation}
acting on $\Omega$ of all degrees. This is easy to see by induction and the (graded)-derivation properties of the ingredients, and establishes the first statement of the lemma for all $\omega\in \Omega$ and $\eta\in \Omega^1$. We now show that if the statement holds for the pairs $(\omega,\eta$) and $(\omega,\eta')$ then it holds for $(\omega,\eta\eta')$ which then establishes the result by induction on the degree of $\eta$. We use again that $\CL_\tau$ is a derivation and now that $\omega\perp$ is a derivation of degree $|\omega|-2$ to break down both sides and compare, using graded-commutativity where needed.These steps are all straightforward and details are again omitted. (2) We next apply  the operations in the first displayed statement to each term of $S(\omega,\eta)$. Since $\extd$ commutes with $\CL_\tau$ we have equality for each term except one, where we have $\alpha\extd(\omega\perp\eta)=\extd(\alpha\omega\perp\eta)- (\extd\alpha)\omega\perp\eta$. The second term results in the right hand side of the second displayed statement. \endproof 

\begin{proposition}\label{tauprop} Let $(\Omega(A),\delta)$ be of classical type and $\tau$ $\delta$-conformal. The following are equivalent

1. $\tau$ is strongly $\delta$-conformal.

2. 
\begin{eqnarray}\label{lieleib}
\CL_\tau L_\delta(\omega,\eta)+\alpha L_\delta(\omega,\eta)-L_\delta(\CL_\tau\omega,\eta)&-&L_\delta(\omega,\CL_\tau\eta)\nonumber\\
&&=-(-1)^{|\omega|}|\omega|\omega\inter_{\extd\alpha}\eta-|\eta|(\inter_{\extd\alpha}\omega)\eta \end{eqnarray}
for all $\omega,\eta\in \Omega$. In this case
\begin{equation}\label{lielap}
[\Delta,\CL_\tau]=\alpha\Delta+(D-\beta)\CL_{\extd\alpha}+(\extd\alpha)\delta +\inter_{\extd\alpha}\extd\end{equation}
where  $\Delta=\extd\delta+\delta\extd$. Moreover, these identities all apply if $(\ ,\ )$ is nondegenerate. \end{proposition}
\proof  (1) We let $C(\omega,\eta)$ first denote the LHS of (\ref{lieleib}). Then
\begin{eqnarray*} \inter_\zeta C(\omega,\eta)&=&\inter_\zeta\CL_\tau L_\delta(\omega,\eta)+\alpha\inter_\zeta L_\delta(\omega,\eta)-\inter_\zeta L_\delta(\CL_\tau\omega,\eta)-\inter_\zeta L_\delta(\omega,\CL_\tau\eta)\\
&=&(\CL_\tau + 2\alpha) \inter_\zeta L_\delta(\omega,\eta)-\inter_{(\CL_\tau\zeta)} L_\delta(\omega,\eta)-\inter_\zeta L_\delta(\CL_\tau\omega,\eta) -\inter_\zeta L_\delta(\omega,\CL_\tau\eta)\\
&=&-(\CL_\tau+2\alpha)\left(L_\delta(\inter_\zeta\omega,\eta)+(-1)^{|\omega|}L_\delta(\omega,\inter_{\zeta}\eta)-L_{\inter_{\extd\zeta}}(\omega,\eta)\right)\\
&&+L_\delta(\inter_{(\CL_\tau\zeta)}\omega,\eta)+(-1)^{|\omega|}L_\delta(\omega,\inter_{(\CL_\tau\zeta)}\eta)-L_{\inter_{\extd\CL_\tau\zeta}}(\omega,\eta)\\
&&+L_\delta(\inter_\zeta\CL_\tau\omega,\eta)+(-1)^{|\omega|}L_\delta(\CL_\tau\omega,\inter_\zeta\eta)-L_{\inter_{\extd\zeta}}(\CL_\tau\omega,\eta)\\
&&+L_\delta(\inter_\zeta\omega,\CL_\tau\eta)+(-1)^{|\omega|}L_\delta(\omega,\inter_\zeta\CL_\tau\eta)-L_{\inter_{\extd\zeta}}(\omega,\CL_\tau\eta)\\
&=&-(\CL_\tau+2\alpha)\left(L_\delta(\inter_\zeta\omega,\eta)+(-1)^{|\omega|}L_\delta(\omega,\inter_{\zeta}\eta)\right)\\
&&+L_\delta(\CL_\tau\inter_\zeta\omega,\eta)+L_\delta(\alpha\inter_\zeta\omega,\eta)+(-1)^{|\omega|}L_\delta(\CL_\tau\omega,\inter_\zeta\eta)\\
&&+L_\delta(\inter_\zeta\omega,\CL_\tau\eta)+(-1)^{|\omega|}L_\delta(\omega,\CL_\tau\inter_\zeta\eta)+(-1)^{|\omega|}L_\delta(\omega,\alpha\inter_\zeta\eta)\\
&=&-(\CL_\tau+\alpha)\left(L_\delta(\inter_\zeta\omega,\eta)+(-1)^{|\omega|}L_\delta(\omega,\inter_{\zeta}\eta)\right)\\
&&+L_\delta(\CL_\tau\inter_\zeta\omega,\eta)+(-1)^{|\omega|-1}(\inter_\zeta\omega) \inter_{\extd\alpha}\eta+(-1)^{|\omega|}L_\delta(\CL_\tau\omega,\inter_\zeta\eta)\\
&&+L_\delta(\inter_\zeta\omega,\CL_\tau\eta)+(-1)^{|\omega|}L_\delta(\omega,\CL_\tau\inter_\zeta\eta)+(-1)^{|\omega|}(\inter_{\extd\alpha}\omega)\inter_\zeta\eta\\
&=&-C(\inter_\zeta\omega,\eta)-(-1)^{|\omega|}C(\omega,\inter_{\zeta}\eta)   -(-1)^{|\omega|}(\inter_\zeta\omega) \inter_{\extd\alpha}\eta  +(-1)^{|\omega|}(\inter_{\extd\alpha}\omega)\inter_\zeta\eta
\end{eqnarray*}
where we used (\ref{lieint}) to swap the derivative with the interior product and  used the first part of Lemma~\ref{Schouten} to break down $L_\delta$. We used 
\[ [\inter_\eta,\CL_\tau]=2\alpha\inter_\eta-\inter_{\CL_\tau\eta},\quad \forall\eta\in \Omega^2,\]
easily deduced from (\ref{lieint}), to cancel the $L_{\inter_{\extd\zeta}}$ and $L_{\inter_{\CL_\tau\extd\zeta}}$ terms arising from Lemma~\ref{Schouten} (where $\CL_\tau$ commutes with $\extd$). Finally, we used Lemma~\ref{deriv} to move $\alpha$ out from inside $L_\delta$ leaving a residue as shown.  We have thus obtained an `induction formula' for $C(\omega,\eta)$ from smaller degrees and this is solved by the right hand side of (\ref{lieleib}). The initial values are $C(a,\omega)=C(\omega,a)=0$ for $a\in A$ and if $(\ ,\ )$ is nondegenerate then the induction step completely determines $C(\omega,\eta)$ as equal to the stated RHS of (\ref{lieleib}). Hence (\ref{lieleib}) holds in the non-degenerate case. (2) Next, let $D$ be the degree operator and observe that the Leibnizator of $\inter_{\extd\zeta}D$ is
\begin{eqnarray*} &&L_{\inter_{\extd\alpha}D}(\omega,\eta)=\inter_{\extd\alpha}D(\omega\eta)-(\inter_{\extd\alpha}D\omega)\eta-(-1)^{|\omega|}\omega \inter_{\extd\alpha}D\eta\\
&=&(|\omega|+|\eta|)((\inter_{\extd\alpha}\omega)\eta+(-1)^{|\omega|}\omega\inter_{\extd\alpha}\eta)-|\omega|(\inter_{\extd\alpha}\omega)\eta-(-1)^{|\omega|}|\eta|\omega \inter_{\extd\alpha}\eta\\
&=&(-1)^{|\omega|}|\omega|\omega\inter_{\extd\alpha}\eta+|\eta|(\inter_{\extd\alpha}\omega)\eta=-C_\alpha(\omega,\eta)
\end{eqnarray*}
where now $C_\alpha(\omega,\eta)$ refers to the expression on the RHS of (\ref{lieleib}). The operator $B=\inter_{\extd\alpha}(D-\beta)$ will have the same Leibnizator since it differs by a constant multiple of $\inter$ and the latter is a derivation, so $L_B(\omega,\eta)=-C_\alpha(\omega,\eta)$. Now suppose that the condition in Definition~\ref{liedelta} holds on $\omega,\eta$, then we show that it holds on $\omega\eta$ provided (\ref{lieleib}) holds:
\begin{eqnarray*} 
&&\delta\CL_\tau(\omega\eta)=\delta((\CL_\tau\omega)\eta)+\delta(\omega\CL_\tau\eta)\\
&=&L_\delta(\CL_\tau\omega,\eta)+L_\delta(\omega,\CL_\tau\eta)+(\delta\CL_\tau\omega)\eta+(-1)^{|\omega|}(\CL_\tau\omega)\delta\eta+(\delta\omega)\CL_\tau\eta+(-1)^{|\omega|}\omega\delta\CL_\tau\eta\\
&=&L_\delta(\CL_\tau\omega,\eta)+L_\delta(\omega,\CL_\tau\eta)+(\CL_\tau\delta\omega)\eta+(-1)^{|\omega|}(\CL_\tau\omega)\delta\eta+(\delta\omega)\CL_\tau\eta+(-1)^{|\omega|}\omega\CL_\tau\delta\eta\\
&&+\alpha(\delta\omega)\eta+(B\omega)\eta+(-1)^{|\omega|}\alpha\omega\delta\eta+(-1)^{|\omega|}\omega B\eta\\
&=&L_\delta(\CL_\tau\omega,\eta)+L_\delta(\omega,\CL_\tau\eta)+(\CL_\tau\delta\omega)\eta+(-1)^{|\omega|}(\CL_\tau\omega)\delta\eta+(\delta\omega)\CL_\tau\eta+(-1)^{|\omega|}\omega\CL_\tau\delta\eta\\
&&+\alpha\delta(\omega\eta)-\alpha L_\delta(\omega,\eta)+B(\omega\eta)+C_\alpha(\omega,\eta)\\
&=&\CL_\tau((\delta\omega)\eta)+(-1)^{|\omega|}\CL_\tau(\omega\delta\eta)+\alpha\delta(\omega\eta)+B(\omega\eta)+\CL_\tau L_\delta(\omega,\eta)\\
&=&\CL_\tau\delta(\omega\eta)+\alpha\delta(\omega\eta)+B(\omega\eta)\\
\end{eqnarray*}
as required. We used the derivation property of $\CL$ and the `compensated derivation' property of $\delta, B$ where we correct with the respective Leibnizator. For the 3rd equality we swapped the order of $\delta,\CL_\tau$ under the inductive assumption of the condition in Definition~\ref{liedelta} and for the 5th equality we use (\ref{lieleib}) and the derivation property of $\CL_\tau$ in reverse. Hence (\ref{lieleib}) implies the condition in Definition~\ref{liedelta} in all degrees. Clearly this proof can be reversed, i.e. if the condition in Definition~\ref{liedelta} holds on all degrees then so does (\ref{lieleib}). (3) Finally, 
\begin{eqnarray*}[\Delta,\CL_\tau]&=&[\delta\extd+\extd\delta,\CL_\tau]=[\delta,\CL_\tau]\extd+\extd[\delta,\CL_\tau]\\
&=&\alpha\delta\extd +\inter_{\extd\alpha}(D-\beta)\extd+\extd(\alpha\delta(\ ))+\extd\inter_{\extd\alpha}(D-\beta)\\
&=&\alpha\Delta-\alpha\extd\delta+\extd(\alpha\delta(\ ))+(D-\beta)\CL_{\extd\alpha}+\inter_{\extd\alpha}\extd
\end{eqnarray*} which we recognise as stated in (\ref{lielap}). 
 \endproof
 
To be sure of all these identities we will require $\tau$ to be strongly $\delta$-conformal. But we see that if $(\ ,\ )$ is nondegenerate then $\delta$-conformal implies strongly $\delta$-conformal. 

\goodbreak
\subsection{Quantisation by a $\delta$-conformal 1-form}
We now combine the two preceding subsections. When $(\Omega(A),\delta)$ is of classical type we have a canonical extension to $\tilde\Omega(A)$ by Theorem~\ref{cleftclasslevi} and a canonical further extension of that to $\tilde{\tilde{\Omega}}(A)$ by Proposition~\ref{Omegafull}. 

\begin{proposition}\label{spacetime} Let $(\Omega(A),\delta)$ be of classical type and  $\tau\in \Omega^1$ strongly $\delta$-conformal. Then $\tau$ defines a  derivation on $\tilde{\tilde{\Omega}}(A)$ in Proposition~\ref{Omegafull} by
\[ \tau(\theta')=\alpha\theta',\quad \tau(\extd\theta')=(\extd\alpha)\theta'+\alpha\extd\theta',\]
\[\tau(\omega)=\CL_\tau\omega +{\lambda\over 2}(-1)^{|\omega|} (|\omega|-\beta)\inter_{\extd\alpha}\omega\theta',\quad\forall \omega\in \Omega\]
The semidirect product DGA $\tilde{\tilde{\Omega}}\lcross\Omega(t,\extd t)$ by Proposition~\ref{semicalc} in the case of a Riemannian manifold essentially recovers  the calculus in \cite[Thm~3.1]{Ma:bh}. 
\end{proposition}
\proof We verify that $\tau$ as defined is a derivation for the $\cdot$ product in Proposition~\ref{Omegafull} that commutes with $\extd_\cdot$. We have already established all needed identities in Proposition~\ref{tauprop}. Thus
\begin{eqnarray*} \tau(\extd_\cdot\omega)&=&\tau\left(\extd \omega - {\lambda\over 2}(-1)^{|\omega|}(\Delta\omega)\theta'+{\lambda\over 2}(\delta\omega)\extd\theta'\right)\\
&=&\CL_\tau\extd\omega +{\lambda\over 2}(-1)^{|\omega|+1} (|\omega|+1-\beta)\inter_{\extd\alpha}\extd\omega\theta'- {\lambda\over 2}(-1)^{|\omega|}\CL_\tau\Delta\omega\theta'+{\lambda\over 2}\CL_\tau\delta\omega\extd\theta'\\
&&- {\lambda\over 2}(-1)^{|\omega|}(\Delta\omega)\alpha\theta'+{\lambda\over 2}(\delta\omega)\extd\alpha\theta'+{\lambda\over 2}(\delta\omega)\alpha\extd\theta'\\
&=&\extd\CL_\tau\omega + {\lambda\over 2}(-1)^{|\omega|+1}(|\omega-\beta)\inter_{\extd\alpha}\extd\omega\theta'- {\lambda\over 2}(-1)^{|\omega|}\Delta\CL_\tau\omega\theta'+{\lambda\over 2}\delta\CL_\tau\omega\extd\theta'\\
&&+ {\lambda\over 2}(-1)^{|\omega|}(|\omega|-\beta)\CL_{\extd\alpha}\omega\theta'-{\lambda\over 2}(|\omega|-\beta)\inter_{\extd\alpha}\omega\extd\theta'\\
&=&\extd\CL_\tau\omega- {\lambda\over 2}(-1)^{|\omega|}\Delta\CL_\tau\omega\theta'+{\lambda\over 2}\delta\CL_\tau\omega\extd\theta'\\
&& +{\lambda\over 2}(-1)^{|\omega|} (|\omega-\beta)\extd\inter_{\extd\alpha}\omega\theta'-{\lambda\over 2}(|\omega|-\beta)\inter_{\extd\alpha}\omega\extd\theta'\\
&=&\extd_\cdot\left(\CL_\tau\omega + {\lambda\over 2}(-1)^{|\omega|} (|\omega|-\beta)\inter_{\extd\alpha}\omega\theta'\right)=\extd_\cdot\tau(\omega)
\end{eqnarray*}
for all $\omega\in\Omega^1$. We used the definitions then strong-$\delta$-conformality and (\ref{lielap}) to reorder. We used graded-commutativity on $(\delta\omega)\extd\alpha$ and we used the Cartan formula for $\CL_{\extd\alpha}$, before recognising the result. Next, we note that $\lo\omega,\eta\lc=L_\delta(\omega,\eta)+S(\omega,\eta)$ according to Theorem~\ref{cleftclasslevi} and we write $B=\inter_{\extd\alpha}(D-\beta)$.  Then
\begin{eqnarray*}\tau(\omega\cdot\eta)&=&\tau\left(\omega\eta+{\lambda\over 2}(-1)^{|\omega|+|\eta|}\lo\omega,\eta\lc\theta'-{\lambda\over 2}(-1)^{|\omega|}(\omega\perp\eta)\extd\theta' \right)\\
&=&\CL_\tau(\omega\eta) +(-1)^{|\omega|+|\eta|} {\lambda\over 2}B(\omega\eta)\theta'\\
&&+{\lambda\over 2}(-1)^{|\omega|+|\eta|}\CL_\tau\lo\omega,\eta\lc\theta'-{\lambda\over 2}(-1)^{|\omega|}\CL_\tau(\omega\perp\eta)\extd\theta'\\
&&+{\lambda\over 2}(-1)^{|\omega|+|\eta|}\lo\omega,\eta\lc\alpha\theta'-{\lambda\over 2}(-1)^{|\omega|}(\omega\perp\eta)(\extd\alpha\theta'+\alpha\extd\theta')\\
&=&\CL_\tau(\omega\eta) + {\lambda\over 2}(-1)^{|\omega|+|\eta|}B(\omega\eta)\theta'+{\lambda\over 2}(-1)^{|\omega|+|\eta|}(\lo\CL_\tau\omega,\eta\lc+\lo\omega,\CL_\tau\eta\lc)\theta'\\
&&-{\lambda\over 2}(-1)^{|\omega|}((\CL_\tau\omega)\perp\eta+\omega\perp\CL_\tau\eta)\extd\theta'+{\lambda\over 2}(-1)^{|\omega|+|\eta|}C_\alpha(\omega,\eta)\theta'\\
&=&(\CL_\tau\omega)\eta+\omega\CL_\tau\eta +{\lambda\over 2}(-1)^{|\omega|+|\eta|}(\lo\CL_\tau\omega,\eta\lc+\lo\omega,\CL_\tau\eta\lc)\theta'\\
&&-{\lambda\over 2}(-1)^{|\omega|}((\CL_\tau\omega)\perp\eta+\omega\perp\CL_\tau\eta)\extd\theta'+ {\lambda\over 2}(-1)^{|\omega|+|\eta|}B(\omega)\eta\theta'+ {\lambda\over 2}(-1)^{|\eta|}\omega B(\eta)\theta'\\
&=& \left(\CL_\tau\omega +{\lambda\over 2}(-1)^{|\omega|} B(\omega)\theta'\right)\cdot\eta+\omega\cdot\left(\CL_\tau\eta +{\lambda\over 2}(-1)^{|\eta|} B(\eta)\theta'\right) \\
&=&\tau(\omega)\cdot\eta+\omega\cdot\tau(\eta)
\end{eqnarray*}
We use the definitions, then the second statement of Lemma~\ref{lieperp} and (\ref{lieleib}) tell us $\CL_\tau\lo\omega,\eta\lc$ while the first of Lemma~\ref{lieperp} tells us $\CL_\tau(\omega\perp\eta)$. We then use from the proof Proposition~\ref{tauprop} that $L_B(\omega,\eta)=-C_\alpha(\omega,\eta)$ and recognise the result.

We then apply Proposition~\ref{semicalc} to obtain a DGA $\tilde{\tilde{\Omega}}(A)\lcross\Omega(t,\extd t)$ over $A_t$ and compute the cross relations for $a\in A$ and $\omega\in \Omega^1$, 
\[ [t,\theta']=\lambda\alpha\theta',\quad[t,\omega]=\lambda(\CL_\tau\omega-\omega)+\lambda^2\left({n-2\over 4}\right)(\extd\alpha,\omega)\theta', \quad [\extd t,a]=\lambda\extd a\]
\[[t,\extd\theta']=\lambda\extd((\alpha-1)\theta'),\quad  \{\extd t,\omega\}=\lambda\extd\omega,\quad \{\extd t,\theta'\}=\lambda\extd\theta',\quad \forall \omega\in\Omega^1\]
which, remembering the change of sign of $\lambda$, is the calculus in \cite[Thm~3.1]{Ma:bh}  and \cite[Prop~3.6]{Ma:bh} in the special  case $\beta=\zeta=0$  in the notation used there and in the case of a Riemannian manifold. \endproof

We have thus put \cite{Ma:bh} into a proper framework and to all degrees: we start with a classical Riemannian manifold and usual exterior algebra $\Omega(M)$, construct  an extension $\tilde{\tilde{\Omega}}(M)$ from Proposition~\ref{Omegafull} and then apply the semidirect product construction Proposition~\ref{semicalc} using a $\delta$-conformal 1-form $\tau$ to give a `extended quantum spacetime' DGA $\tilde{\tilde{\Omega}}(M)\lcross\Omega(t,\extd t)$. The special case covered corresponds to a deformation of a direct product metric on $M\times\R$  whereas the full generality in \cite{Ma:bh} can be expected to require a cocycle semidirect product version of Proposition~\ref{semicalc}. An example given in detail in \cite{Ma:bh} is the Schwarzschild black hole.

\begin{remark} Starting with any standard DGA $\Omega(A)$ it is straightforward to formulate a general version of Proposition~\ref{Omegafull} with $\Delta,\delta,\lo\ ,\ \lc, \perp$ graded maps of the appropriate degree used to define an associated product $\cdot$ and differential $\extd_\cdot$ and to write out the requirements on these maps to arrive at a DGA. This is basically the composition of Theorem~\ref{construct} and Proposition~\ref{Omegafull} and for this reason the details are omitted but it is clear from the form of the commutation relations that Proposition~\ref{spacetime} provides a non-graded-commutative example where the role of $\Omega(A)$ is now played by the standard inner DGA $\Omega(A)\lcross\Omega(t,\extd t)$. This is because we can use the commutation relations to move  $\theta',\extd\theta'$ to the far right. 
\end{remark}

\begin{remark} One can take $\tau=0$ in Proposition~\ref{semicalc}. Then any DGA on $A$ has a canonical extension to $\Omega(A[t])=\Omega(A)\lcross\Omega(t,\extd t)$ as a DGA on the central extension $A[t]$. Here the cross relations are
\[ [\omega,t]=\lambda|\omega|\omega,\quad [\extd t,\omega\}=\lambda\extd \omega,\quad \forall \omega\in\Omega(A).\]
One can check that these relations also hold for $\omega\in \Omega(A[t])$, i.e. the canonical degree derivation $D$ and exterior graded-derivation $\extd$ on $\Omega(A)$ are rendered inner by $t$ and $\extd t$ respectively. 

For example, iterating this case  from $k$ we have $\Omega(t_1,\extd t_1)\lcross \Omega(t_2,\extd t_2)\cdots\lcross \Omega(t_n,\extd t_n)$ as a quantum calculus on $k[t_1,\cdots,t_n]$, which can be written as
\[ \Omega=k[t_1,\cdots,t_n]\lcross\Lambda(\extd t_1,\cdots,\extd t_n),\quad  [\extd t_i,t_j]=\lambda\extd t_{{\rm min}(i,j)}\]
 where $\Lambda$ is the usual Grassmann algebra on $n$ variables $\extd t_i$. By construction, this is inner with $\theta=\extd t_n$ and has inner grading by commutator with $t_n$.
 
Similarly for any $M$ manifold by taking $\tau=0$ in Proposition~\ref{semicalc} we have a non-graded-commutative DGA,  $\Omega(M)\lcross\Omega(t,\extd t)$ on  $C(M)[t]$, the latter being a classical  commutative algebra of functions on $M\times\R$.   \end{remark}


\end{document}